\documentclass[final,3p,times]{elsarticle}
\usepackage{amssymb}
\usepackage{amsmath}
\usepackage{amsthm}

\usepackage{amsgen,amstext}
\usepackage{amsfonts,bm}
\usepackage{dsfont}
\usepackage{mathrsfs}
\usepackage{latexsym}
\usepackage[utf8]{inputenc} 
\usepackage[T1]{fontenc}    
\usepackage{hyperref}       
\usepackage{xurl}
\usepackage{booktabs}       
\usepackage{nicefrac}       
\usepackage{microtype}      
\usepackage{caption,subcaption}
\usepackage{xcolor}
\usepackage[section]{algorithm}
\usepackage{algcompatible}

\newtheorem{theorem}{Theorem}[section]
\newtheorem{lemma}[theorem]{Lemma}

\DeclareMathOperator*{\diag}{diag}

\newcommand{\Ints}{\mathbb{Z}}
\newcommand{\Nats}{\mathbb{N}}
\newcommand{\Reals}{\mathbb{R}}
\newcommand{\pr}{\mathbb{P}}

\newcommand{\bb}{\mathbf{b}}
\newcommand{\bc}{\mathbf{c}}
\newcommand{\bd}{\mathbf{d}}

\newcommand{\bu}{\mathbf{v}}
\newcommand{\bw}{\mathbf{s}}
\newcommand{\bx}{\mathbf{x}}
\newcommand{\byy}{\mathbf{y}}
\newcommand{\by}{\bm{\pi}}
\newcommand{\bz}{\mathbf{z}}
\newcommand{\bA}{\mathbf{A}}
\newcommand{\bD}{\mathbf{D}}

\newcommand{\bH}{\mathbf{T}}
\newcommand{\bI}{\mathbf{I}}
\newcommand{\bJ}{\mathbf{B}}
\newcommand{\bL}{\mathbf{L}}
\newcommand{\bM}{\mathbf{M}}
\newcommand{\bP}{\mathbf{P}}
\newcommand{\bQ}{\mathbf{R}}
\newcommand{\bT}{\mathbf{W}}
\newcommand{\bU}{\mathbf{V}}
\newcommand{\bW}{\mathbf{S}}
\newcommand{\bX}{\mathbf{X}}
\newcommand{\bZ}{\mathbf{Y}}
\newcommand{\bzero}{\mathbf{0}}
\newcommand{\ones}{\mathbf{1}}

\newcommand{\bUU}{\mathbf{U}}

\newcommand{\eb}{{b}}
\newcommand{\ec}{{c}}
\newcommand{\eh}{{w}}
\newcommand{\es}{{h}}
\newcommand{\et}{{g}}
\newcommand{\ew}{{s}}

\begin{document}

\begin{frontmatter}

\title{On Mixed-Precision Iterative Methods and Analysis \\for Nearly Completely Decomposable Markov Processes}

\author{Vasileios Kalantzis, Mark S.\ Squillante, Chai Wah Wu} 

\affiliation{
organization={Mathematics of Computation, IBM Research},
            addressline={Thomas J.\ Watson Research Center}, 
            city={Yorktown Heights},
            postcode={10598}, 
            state={NY},
            country={USA}}

\begin{abstract}
In this paper we consider the problem of computing the stationary distribution of nearly completely decomposable Markov processes, a well-established area in the classical theory of Markov processes with broad applications in the design, modeling, analysis and optimization of computer systems and applications.
We
devise a general mathematical framework of numerical
solution 
methods
that exploits forms of mixed-precision computation to significantly reduce computation times
and that exploits forms of iterative approximate 
computing approaches
to mitigate the 
impact of inaccurate computations, further reduce computation times, and ensure convergence.
Then we derive a mathematical analysis 
that establishes theoretical 
properties of our general algorithmic framework including results
on approximation errors, convergence behaviors, and other algorithmic 
characteristics.
Numerical experiments 
demonstrate that our general algorithmic 
framework provides
significant improvements in computation times over the most-efficient existing numerical methods.
\end{abstract}


\end{frontmatter}

\section{Introduction}
Markov processes often play an important role in the design and mathematical modeling, analysis and optimization of a wide range of computer systems
and applications.
A critically important goal in the mathematical modeling, analysis and optimization of these computer systems and
applications as well as
their dynamic behavior concerns the need to determine the stationary distribution of the corresponding Markov process, from which various measures and quantities of interest can 
in turn
be directly obtained to facilitate and support such design, analysis and optimization of 
these
computer systems and applications.

More precisely, consider a discrete-time Markov process $\{ X(s) \, ; \, s \in \Ints_+ \}$ defined on the state space $[n] := \{ 1, \ldots, n \}$,
and let ${\bP} \in \Reals^{n\times n}$ and $\by := (\pi_1, \ldots, \pi_n)$ respectively denote its transition probability matrix and stationary distribution.
Here, $\Reals$ and $\Ints_+$ denote the set of real numbers and nonnegative integers, respectively.
Assuming this Markov process is irreducible and ergodic, then the invariant probability vector $\by$ exists and is uniquely determined as the solution of~\cite{Stewart94,Stroock14}
\begin{equation}
    \by = \by \bP \qquad \mbox{and} \qquad \by \ones = \| \by \|_1 = 1 ,
\label{eq:2.4}
\end{equation}
where 
$\ones = (1,\ldots,1)^\top$ 
is the column vector of appropriate dimension containing all ones,
and $\| \cdot \|_1$ is the $l_1$ norm of any given vector.
While directly computing the solution of~\eqref{eq:2.4} is feasible for moderate to relatively large values of $n$, the state space size $n$ is often very large for the Markov processes of many current and emerging computer systems and applications, thus causing the direct computation of the solution of~\eqref{eq:2.4} to become prohibitively expensive.
In these frequently encountered situations, one needs to exploit structural properties of the Markov process and its underlying computer systems and applications to realize more efficient computation.

One particularly important class of structural properties often arising in such Markov processes concerns the so-called ``aggregation of variables''~\cite{SimAnd61} in which all the states (variables) of the system can be grouped into a relatively small number of relatively large clusters such that the magnitude of 
transition probabilities
between states within a cluster overwhelmingly dominates the relatively tiny magnitude of 
transition probabilities
between states in different clusters.
Markov processes with these structural properties, called \emph{nearly completely decomposable} (NCD), make it possible to effectively investigate the system and its dynamics by exploiting the analysis of interactions among the states within each individual cluster without reference to the interactions among the states across clusters and by exploiting the analysis of interactions among the states across different clusters without reference to the interactions among the states within clusters.
Such NCD properties arise naturally in many computer systems and applications due to the prevalence of levels of modular abstraction or hierarchies of components in the design of these computer systems and applications, where interactions within modules or components dominate interactions between modules or components.
In addition to the various traditional Markov process instances within computer systems
and applications
(see, e.g.,~\cite{Cour75,CouVan76,Cour77,AvCoKo87}), these NCD structural properties arise naturally in various traditional aspects of many other domains such as
queueing theory,
control theory,
dynamical systems,
inventory and production systems,
economics, 
biology,
genetics,
and
social sciences; refer to, e.g.,~\cite{ArKaSc58,SimAnd61,FisAnd62,AndFis63,Cour77,AvCoKo87,AldKha91,YinZha05}.

More formally, consider
a
discrete-time Markov process $\{ X(s) \, ; \, s \in \Ints_+ \}$ 
defined on the state space $[n]$
with transition probability matrix $\bP$ taking the block structural form
\begin{equation}
    \bP = \begin{pmatrix}
        \bP_{11} & \bP_{12} & \cdots & \bP_{1m} \\
        \bP_{21} & \bP_{22} & \cdots & \bP_{2m} \\
        \vdots & \vdots & \ddots & \vdots \\
        \bP_{m1} & \bP_{m2} & \cdots & \bP_{mm}
        \end{pmatrix} ,
\label{eq:P-matrix}
\end{equation}
where
$\bP_{ij}$, $i,j\in[m]$, are nonnegative block matrices in $\Reals^{n_i \times n_j}$ such that $\sum_{j\in[m]} \bP_{ij}$ are stochastic
(i.e., rows sum to one with each element in $[0,1]$),
and 
$\bP$ is a stochastic matrix in $\Reals^{n \times n}$, $n = \sum_{i\in[m]} n_{i}$.
Define $\sigma_i := \sum_{\ell=1}^{i-1} n_{\ell}$, 
$i \in [m]$.
Under the assumption that the Markov process is NCD, the elements of the off-diagonal block matrices $\bP_{ij}$ are very small relative to the elements of the diagonal block matrices $\bP_{ii}$,
$i,j \in [m], i\neq j$.
In particular, we assume throughout that
\begin{equation}
    \| \bP_{ii} \| = O(1) \qquad \mbox{ and } \qquad \| \bP_{ij} \| = O(\epsilon) , \qquad\qquad i,j \in [m], \; i\neq j ,
\label{eq:2.2}
\end{equation}
where $\| \cdot \|$ is the spectral norm of any given matrix and $\epsilon > 0$ is a sufficiently small constant.
Define
\begin{equation}
\by := (\by_1, \by_2, \ldots, \by_m) ,
\qquad
\by_i := (\pi_{i1}, \pi_{i2}, \ldots, \pi_{i n_i}),
\qquad
\pi_{ik} := \lim_{s \rightarrow \infty} \pr[ X(s) = (\sigma_i+k) ] ,
\label{eq:pi-def}
\end{equation}
for $k \in [n_i]$, $i \in [m]$.
The limiting probability vector $\by$ is the stationary distribution of the Markov process $\{ X(s) \, ; \, s \in \Ints_+ \}$.
We assume throughout that this process is irreducible and ergodic, and thus its invariant probability vector $\by$ exists and is uniquely determined as the solution of~\eqref{eq:2.4}.

In addition to their broad applications across a wide range of domains in computing, engineering, science and business, NCD structural properties represent an important fundamental area in the classical theory of Markov processes~\cite{SimAnd61,Taka75,Vant81,KoMcSt84,McStSt84,CaoSte85,Havi86,StStMc94,YinZha05}.
Hence, several different (yet related) numerical methods have been developed to efficiently compute the stationary distribution $\by$ of NCD Markov processes, with the most efficient based on different forms of exploiting the aggregation and disaggregation of NCD structural properties as part of an iterative computational method.
This collection of the most-efficient numerical methods for NCD Markov processes consists predominantly of the algorithms due to Takahashi~\cite{Taka75}, Vantilborgh~\cite{Vant81}, and Koury, McAllister and Stewart~\cite{KoMcSt84}.
On the one hand, although they differ in terms of specific details, 
all
these methods exploit aggregation-disaggregation as part of an iterative computational process in a very similar manner and with similar convergence behaviors~\cite{CaoSte85}.
On the other hand, despite their significant performance benefits over directly solving~\eqref{eq:2.4}, these numerical methods can still be prohibitively expensive for computing the stationary distribution of highly large-scale NCD Markov processes and for supporting online or real-time applications of such stochastic processes.
Moreover, all of these most-efficient methods involve computations in full precision and none of them take advantage of 
advances in
mixed-precision technology.

\textbf{Our Contributions.}
Our focus in this paper is on the design and analysis of 
a general mathematical framework of
numerical methods for computing the stationary distribution $\by$ of NCD Markov processes that addresses the performance bottlenecks of existing 
algorithms
and provides significant computational improvements to support the very large scale of current and emerging computer systems and applications as well as their potential deployment in online or real-time environments.
In doing so, we make the following important technical contributions.
We first devise 
a
general 
algorithmic framework
of computational 
approaches to efficiently determine the solution of~\eqref{eq:2.4} for the invariant probability vector $\by$ in~\eqref{eq:pi-def} of the Markov process $\{ X(s) \, ; \, s \in \Ints_+ \}$ whose $\bP$ matrix has the structural form~\eqref{eq:P-matrix} with properties~\eqref{eq:2.2}.
Our general mathematical framework
involves a combination of different forms of exploiting advances in computer architectures related to \emph{mixed-precision computation},
which significantly reduce computation times at the expense of inaccuracies,
and exploiting 
advances in
\emph{iterative approximate computing methods},
which mitigate the impact of inaccurate computations, further reduce computation times, and guarantee convergence.
It is important to note that, although mixed-precision computation has received considerable attention in machine learning and related applications, such approaches and methods are far less well established in the context of numerical solutions of stochastic processes;
moreover, our framework combines such mixed-precision computation with iterative approximate computing methods.
We then derive a mathematical analysis of our general 
algorithmic
framework
that rigorously establishes their important theoretical properties including results on 
approximation errors, convergence behaviors, and other algorithmic characteristics.
Lastly, we conduct many numerical experiments which empirically demonstrate that our general algorithmic 
framework provides
orders of magnitude improvements in the computation times to obtain the stationary distribution of NCD Markov processes over the most-efficient existing methods which exploit aggregation-disaggregation.
These numerical experiments are based on simulation of different graphics processing unit (GPU) architectures that we validated against physical experiments on a particular GPU architecture.

Meanwhile, we note that Horton and Leutenegger~\cite{HorLeu94}, motivated by multigrid methods, propose an algorithmic approach for computing the stationary distribution of Markov processes in general without any special structural properties.
Although the proposed algorithm is somewhat related to aggregation-disaggregation algorithms, with the latter methods seen as a special case under specific algorithmic choices, Horton and Leutenegger do not provide any theoretical results on convergence or other algorithmic properties, and they only provide empirical results for Markov processes without NCD structural properties~\cite{HorLeu94}.
For these reasons, we do not consider further in this paper the multigrid-based algorithm of Horton and Leutenegger.
It is important to note, however, that
the numerical methods of our general
mathematical framework
devised in this paper within the context of aggregation-disaggregation algorithms can be analogously exploited within the context of their multigrid-based algorithm for more general Markov processes beyond those having NCD structural properties.

\textbf{Organization of paper.}
Section~\ref{sec:prelim} provides technical 
preliminaries on NCD Markov processes, followed by the presentation of our design of
a general mathematical
solution 
framework together with a few of its representative examples
in Section~\ref{sec:algs}.
We derive our mathematical analysis of these
examples of our
algorithmic
framework
in Section~\ref{sec:analysis}, followed by our proofs of the corresponding theoretical results in Section~\ref{sec:proofs}.
The empirical results from a representative sample of 
our
numerical experiments are presented in Section~\ref{sec:experiments}, followed by concluding remarks.

\section{Technical Preliminaries}
\label{sec:prelim}
In this section, we present preliminary technical details and results for NCD Markov processes, including additional standard conditions and properties assumed throughout the paper.
We shall use the following standard notation.
Let $\Reals$, $\Ints_+$ and $\Nats$ denote the set of real numbers, nonnegative integers and natural numbers, respectively.
Define $[n] := \{1,\ldots,n\} = \Nats \setminus \{ n+1, n+2, \ldots \}$.
Bold lowercase variables denote vectors, while bold uppercase variables denote matrices.
Let $\bI$ denote the identity matrix of appropriate dimension with $\bI_n$ denoting the $n\times n$ identity matrix.

Beyond our assumptions that the Markov process $\{ X(s) \, ; \, s \in \Ints_+ \}$ has the structural form~\eqref{eq:P-matrix} with properties~\eqref{eq:2.2} and is irreducible and ergodic, we further assume that $\tau \leq \|\by_i\|_1$ for some constant $\tau>0$, and thus from~\eqref{eq:2.4} we have $0 < \tau \leq \|\by_i\|_1 < 1$, $i\in [m]$.
This assumption of $\tau \leq \|\by_i\|_1$ simply ensures that none of the $\by_i$ go to $0$ as $\epsilon$ goes to $0$, $i\in [m]$, which is a standard regularity condition~\cite{StStMc94}.
The transition probability matrix $\bP$, since it is a stochastic matrix, has a dominant unit eigenvalue and corresponding right eigenvector $\ones$ given by
$\bP \ones = \ones$.
We additionally assume no other eigenvalues of $\bP$ have absolute value one.

Turning to the properties of the eigenvalues of the diagonal block matrices, let $\bu_i > 0$ denote the left eigenvector of $\bP_{ii}$ that corresponds to its dominant eigenvalue $\lambda_{i_1}$, $i\in[m]$.
We assume that
\begin{equation*}
    \lambda_{i_1} = 1 - \et_i \epsilon + o(\epsilon) , \qquad \et_i > 0, \;\; i \in [m],
\end{equation*}
i.e., 
the dominant eigenvalue of $\bP_{ii}$ is close to one.
Let us further
assume that the remaining eigenvalues of $\bP_{ii}$ are uniformly bounded away from one, whose underlying condition guaranteeing this assumption will be justified below.
Then, it is well known~\cite{McStSt84} that 
there exists a nonsingular matrix $\bT_i = \begin{pmatrix} \bu_i \\ \bU_i \end{pmatrix}$
such that
\begin{equation}
    \begin{pmatrix} \bu_i \\ \bU_i \end{pmatrix} \bP_{ii} = \begin{pmatrix} 1 - \et_i \epsilon + o(\epsilon) & 0 \\ \bzero & \bH_i \end{pmatrix} \begin{pmatrix} \bu_i \\ \bU_i \end{pmatrix} ,
    \label{eq:4.1}
\end{equation}
with the norms of $\bT_i$ and $\bT_i^{-1}$ bounded by a constant independent of $\epsilon$, $i\in [m]$.
From~\eqref{eq:4.1}, we derive
\begin{equation}
\begin{pmatrix} \bu_i \\ \bU_i \end{pmatrix} (\bI - \bP_{ii}) = \begin{pmatrix} \et_i \epsilon + o(\epsilon) & 0 \\ \bzero & \bI - \bH_i \end{pmatrix} \begin{pmatrix} \bu_i \\ \bU_i \end{pmatrix} 
\qquad \mbox{ and } \qquad
\begin{pmatrix} \bu_i \\ \bU_i \end{pmatrix} (\bI - \bP_{ii})^{-1} = \begin{pmatrix} (\et_i \epsilon + o(\epsilon))^{-1} & 0 \\ \bzero & (\bI - \bH_i)^{-1} \end{pmatrix} \begin{pmatrix} \bu_i \\ \bU_i \end{pmatrix} . \label{eq:4.1:cond3}
\end{equation}
By
imposing the condition
\begin{equation}
\| (\bI - \bH_i)^{-1} \| \leq \eh_i ,
\label{eq:4.2}
\end{equation}
for a constant $\eh_i$ independent of $\epsilon$,
we guarantee that the nondominant eigenvalues $\lambda_{i_2}, \ldots, \lambda_{i_{n_i}}$ of $\bP_{ii}$ are uniformly bounded away from one, $i\in [m]$, as assumed above.

Now, define
\begin{equation}
\hat{\by} := (\hat{\by}_1, \hat{\by}_2, \ldots, \hat{\by}_m ) ,
\qquad
\hat{\by}_i := \frac{\by_i}{\| \by_i \|_1}, 
\qquad
\bQ_{ij} := \hat{\by}_i \bP_{ij} \ones ,
\label{eq:2.6}
\end{equation}
for $i,j \in [m]$.
It follows from our assumptions on the matrix $\bP$ and from $\| \hat{\by}_i \|_1 = 1$, $i \in [m]$, that $\bQ$ is an irreducible stochastic matrix with unique left eigenvector $\bw > 0$ corresponding to the unit eigenvalue.
More precisely, we have
\begin{equation}
    \bw = \bw \bQ
    \qquad \mbox{and} \qquad
    \bw \ones = \| \bw \|_1 = 1 .
    \label{eq:2.8}
\end{equation}
Further observe that
\begin{equation*}
\| \by_j \|_1 = \by_j\ones = \sum_{i=1}^m \by_i \bP_{ij} \ones = \sum_{i=1}^m \| \by_i \|_1 \frac{\by_i}{\| \by_i \|_1} \bP_{ij} \ones = \sum_{i=1}^m \| \by_i \|_1 \bQ_{ij}
\end{equation*}
together with $\sum_{i=1}^m \| \by_i \|_1 = 1$, and thus the vector $(\| \by_1 \|_1 , \| \by_2 \|_1 , \ldots, \| \by_m \|_1 )$ satisfies~\eqref{eq:2.8}, which as a result of the uniqueness of its solution yields
\begin{equation}
    \bw = (\| \by_1 \|_1 , \| \by_2 \|_1 , \ldots, \| \by_m \|_1 ) .
    \label{eq:2.9}
\end{equation}

Next, it follows from $\bP$ being an irreducible matrix that there exist uniformly bounded matrices $\bX$ and $\bZ$ such that
\begin{equation}
    \bW := (\ones , \bX) = \begin{pmatrix} \by \\ \bZ \end{pmatrix}^{-1} ,
    \qquad
    \bW^{-1} \bQ \bW = \begin{pmatrix} 1 & 0 \\ 0 & \bJ \end{pmatrix} ,
    \qquad
    \| \bW \| = O(1) =  \| \bW^{-1} \| .
\label{eq:4.4} 
\end{equation}
By further assuming
\begin{equation}
    \| (\bI - \bJ)^{-1} \| = O\left(\frac{1}{\epsilon}\right) ,
\label{eq:4.5}
\end{equation}
we guarantee the nondominant eigenvalues of the aggregation matrix $\bQ$ are bounded away from one by an amount proportional to $\epsilon$.

\section{Algorithmic Solutions}
\label{sec:algs}
In this section, we 
first devise a general mathematical framework of numerical methods
for computing the solution of~\eqref{eq:2.4} to determine the stationary distribution $\by$ of NCD Markov processes with the goal of providing significant improvements in computational and theoretical properties over the most-efficient numerical methods for doing so.
While our general algorithmic framework
of computational approaches 
can be applied within the context of any of these
existing
numerical methods and beyond, we then 
consider representative examples of applying our algorithmic framework within the context of
the algorithm due to Koury, McAllister and Stewart~\cite{KoMcSt84}
(KMS), 
since it is
the most recent of the
most-efficient numerical methods.
In particular, we consider two such examples of our general algorithmic framework
in which we sharply reduce the computational bottlenecks associated with solving systems of linear equations by exploiting
general classes of iterative approximate mixed-precision computing methods.
The primary differences between 
these two representative examples of our algorithmic framework
concern the 
degree of
aggressiveness
with
which we exploit the
general classes of iterative approximate mixed-precision computing methods.

\subsection{General Mathematical Framework}
All the computational approaches to compute the stationary distribution of NCD Markov processes, including the most-efficient aggregation-disaggregation algorithms such as KMS, involve numerical methods to solve systems of linear equations, generically denoted here by $\bA\bx=\bb$,
where the numerical solution of such 
linear systems represents the core computational bottlenecks of these  approaches.
Our general mathematical framework consists of a combination of universal forms of exploiting: 
(1) advances in computer architectures and technology related to mixed-precision computation,
which significantly reduce computation times at the expense of inaccurate results;
and
(2) advances in iterative approximate computing methods,
which mitigate the impact of inaccurate computations, further reduce computation times, and guarantee convergence.
More precisely, each individual invocation of standard full-precision linear system solvers to compute occurrences of $\bA\bx=\bb$ is replaced by a sequence of computational steps based on general classes of iterative approximate mixed-precision methods for numerically solving linear systems.
This involves each computation of $\bA\bx=\bb$ in a chosen mixed-precision mode supported by advanced computer architectures and technology as part of a selected iterative approximate numerical computing method.
The former chosen mixed-precision mode significantly reduces computation times at the cost of inaccurate results, instances of which include forms of multi-precision arithmetic or stochastic rounding;
whereas the latter selected iterative approximation numerical method addresses the inaccurate computations of mixed-precision modes, further reduces computation times and ensures convergence, instances of which include forms of 
iterative refinement (IR)~\cite{Wilkinson}, 
Richardson iteration (RI)~\cite{golub1961chebyshev}, and 
generalized minimal residual (GMRES)~\cite{gmres} methods.

Since both 
representative examples
of 
applying
our general algorithmic framework herein are conducted within the context
of the 
KMS
algorithm, we first present the main steps the original algorithm.
In particular, the KMS algorithm is composed of an iterative numerical method that starts at iteration~$0$ with any given initial approximation $\by^{(0)} = (\by_1^{(0)}, \ldots, \by_m^{(0)})$ of the desired solution $\by$.
Then, for each subsequent iteration~$t=1,2,3,\ldots$:\\
\textbf{Step}~$\mathbf{1}$ consists of normalizing the vector components of the current estimate $\by^{(t-1)}$ of the solution $\by$, i.e., compute $\hat{\by}_i^{(t-1)} = \by_i^{(t-1)} / \| \by_i^{(t-1)} \|_1$ according to the middle equation in~\eqref{eq:2.6}, $i\in[m]$.\\
\textbf{Step}~$\mathbf{2}$ comprises obtaining the elements of the aggregation matrix $\bQ^{(t-1)}$ leading up to the current iteration, i.e., compute $\bQ_{ij}^{(t-1)} = \hat{\by}_i^{(t-1)} \bP_{ij} \ones$ according to the rightmost equation in~\eqref{eq:2.6}, $i,j\in[m]$.\\
\textbf{Step}~$\mathbf{3}$ consists of obtaining the dominant left eigenvector $\bw$ of $\bQ$ by computing the solution of $\bw^{(t-1)} = \bw^{(t-1)} \bQ^{(t-1)}$ and $\bw^{(t-1)} \ones = \| \bw^{(t-1)} \|_1 = 1$ according to~\eqref{eq:2.8}.\\
\textbf{Step}~$\mathbf{4}$ comprises computing the Hadamard product $\bz^{(t)} = \bw^{(t-1)} \odot \hat{\by}^{(t-1)} = (s_1^{(t-1)} \hat{\by}_1^{(t-1)}, \ldots, s_m^{(t-1)} \hat{\by}_m^{(t-1)})$. \\
\textbf{Step}~$\mathbf{5}$ consists of solving the series of $m$ systems of linear equations $\by_i^{(t)} = \by_i^{(t)} \bP_{ii} + \sum_{j<i} \bz_j^{(t)} \bP_{ji} + \sum_{j>i} \by_j^{(t)} \bP_{ji}$ in the order $i=m,\ldots,1$, thus rendering the blockwise estimates of the components of $\by^{(t)}$ for iteration~$t$.\\
\textbf{Step}~$\mathbf{6}$ comprises
a test for convergence,
halting with the return of the current estimate $\by^{(t)}$ if the solution is sufficiently accurate; otherwise, the iteration index $t$ is incremented and the iterative process of \textbf{Step}~$\mathbf{1}$~--~\textbf{Step}~$\mathbf{6}$ is repeated.

The core performance bottlenecks of the KMS algorithm concern computing the solution of $m+1$ systems of linear equations in each iteration~$t$, specifically the series of $n_i\times n_i$ linear systems in \textbf{Step}~$\mathbf{5}$, $i\in[m]$, as well the $m\times m$ linear system in \textbf{Step}~$\mathbf{3}$ for problems with large $m$.
In strong contrast, \textbf{Step}~$\mathbf{1}$ and \textbf{Step}~$\mathbf{2}$ simply set up the computation in \textbf{Step}~$\mathbf{3}$ and involve relatively little computation  
compared to the solution of such a linear system in \textbf{Step}~$\mathbf{3}$;
similarly, \textbf{Step}~$\mathbf{4}$ simply sets up the computation in \textbf{Step}~$\mathbf{5}$ and involves relatively little computation 
compared to the solution of such linear systems in \textbf{Step}~$\mathbf{5}$.
We now present in turn two representative examples of applying our general algorithmic framework to address the core performance bottlenecks of 
\textbf{Step}~$\mathbf{5}$ and \textbf{Step}~$\mathbf{3}$.

\subsection{Representative Conservative Example of General Algorithmic Framework}
\label{ssec:algs-conservative}
The
first 
representative example of a
general 
application of our algorithmic framework within the context
of the KMS algorithm consists of replacing 
each invocation of a full-precision linear system solver
in \textbf{Step}~$\mathbf{5}$ and \textbf{Step}~$\mathbf{3}$ 
with a combination of the selection of an appropriate level of mixed-precision supported by the available computer architecture and the selection of any appropriate instance of the general classes of iterative approximate numerical methods that jointly realize full-precision results.
This representative example of our general algorithmic framework is conservative in the sense that the iterative approximate mixed-precision computing method solves each occurrence of a linear system through a sequence of iterations which lead to full-precision accuracy (though more efficiently).
We summarize such a representative conservative example in Algorithm~\ref{alg:mkms} where the IR method is chosen for the iterative approximate numerical computing method of our general mathematical framework.

As the initial step comprising mixed-precision mode selection, in order to
reduce the computation costs as much as possible while still realizing full-precision results, we determine a level of reduced precision (such as single-precision, half-precision, bfloat16, and so on) based on the 
numerical
properties of the linear system that ensures convergence under the 
IR method
to yield full-precision accuracy.
For our representative example employing the IR method,
the selection of the proper reduced precision is based on an estimate of the condition number and the norm of the matrix $\bP_{ii}$~\cite{HornJohnson13,Strang23} for each linear system in \textbf{Step}~$\mathbf{5}$, $i\in[m]$, and similarly for the linear system of \textbf{Step}~$\mathbf{3}$.
Specifically, for each linear system 
generically denoted by $\bA\bx=\bb$,
we determine the machine precision $\delta$ such that 
$\delta \kappa(\bA)\|\bA\|$
is sufficiently small to ensure convergence~\cite{Moler}, where 
$\kappa(\bA)$
denotes the condition number of the matrix 
$\bA$.

Once the level of reduced precision has been determined in this manner, the 
IR
method is used to solve each of the linear systems in \textbf{Step}~$\mathbf{5}$ and \textbf{Step}~$\mathbf{3}$, 
again
generically denoted by $\bA\bx=\bb$.
In particular, at every step $k$ of the 
IR
method within each outer-loop iteration~$t$ of Algorithm~\ref{alg:mkms}, we compute the residual $\bd_k=\bb-\bA\bx_k$ and then solve for $\byy_k$ in $\bA \byy_k = \bd_k$ using the 
selected
reduced-precision method.
The resulting solution $\byy_k$ is added to the previous estimate $\bx_k$ to obtain $\bx_{k+1} = \bx_k + \byy_k$, and the step index $k$ is incremented.
This process of inner-loop steps is repeated until the residual falls below the error of double-precision arithmetic, which is guaranteed by our selection of the reduced-precision level~\cite{Moler}.
We note that an
analogous approach can be followed for any of the iterative approximate mixed-precision computing method alternatives to
the IR method.

\begin{algorithm}[htpb!]
\caption{Representative Conservative Example of General Algorithmic Framework} \label{alg:mkms}
\begin{algorithmic}
 \State {\bf Input:} Irreducible stochastic matrix $\bP$ with structural form~\eqref{eq:P-matrix}
  \State {\bf Output:} Estimate of probability distribution $\by$ such that $\by = \by \bP$ and $\by \ones = 1$

  \State Set $t=1$ and Let $\by^{(0)}$ be the given initial approximation to the solution $\by$

\State{\bf While:} Stopping criteria is not satisfied
{
\STATE \quad $\mathbf{1.}$ Compute $\hat{\by}_i^{(t-1)} = \frac{\by_i^{(t-1)}}{\| \by_i^{(t-1)} \|_1}$, $i\in[m]$, and Obtain $\hat{\by}^{(t-1)} = (\hat{\by}_1^{(t-1)}, \ldots, \hat{\by}_m^{(t-1)})$
\State \quad $\mathbf{2.}$ Compute $\bQ_{ij}^{(t-1)} = \hat{\by}_i^{(t-1)} \bP_{ij} \ones$, $i,j\in[m]$, and Obtain $\bQ^{(t-1)}$
\State \quad $\mathbf{3.}$ Compute solution of linear system $\bw^{(t-1)} = \bw^{(t-1)} \bQ^{(t-1)}$ and $\bw^{(t-1)} \ones = \| \bw^{(t-1)} \|_1 = 1$ using approximate mixed- \hspace*{0.4in} precision computing method of interest
\State \quad $\mathbf{4.}$ Compute Hadamard product $\bz^{(t)} = \bw^{(t-1)} \odot \hat{\by}^{(t-1)}$
\State \quad $\mathbf{5.}$ Compute solutions of linear systems $\by_i^{(t)} = \by_i^{(t)} \bP_{ii} + \sum_{j<i} \bz_j^{(t)} \bP_{ji} + \sum_{j>i} \by_j^{(t)} \bP_{ji}$ for $i=m,\ldots,1$ using approximate \hspace*{0.4in} mixed-precision computing method of interest, and Obtain $\by^{(t)}$
\State \quad $\mathbf{6.}$ Conduct convergence test, Setting $t=t+1$ if stopping criteria not satisfied
}
\State{\bf Return:} Estimate $\by^{(t)}$
\end{algorithmic}
\end{algorithm}

\subsection{Representative Aggressive Example of General Algorithmic Framework}
\label{ssec:algs-aggressive}
The second representative example of a general application of our algorithmic framework within the context of the KMS algorithm also consists of replacing each invocation of a full-precision linear system solver in \textbf{Step}~$\mathbf{5}$ and \textbf{Step}~$\mathbf{3}$ with a combination of the selection of an appropriate level of mixed-precision supported by the available computer architecture and the selection of any appropriate instance of the general classes of iterative approximate numerical methods that may not necessarily jointly lead to full-precision accuracy.
Specifically, as a primary difference with Section~\ref{ssec:algs-conservative}, we take a more aggressive approach in this representative example of our general algorithmic framework
by exploiting the general classes of iterative approximate mixed-precision numerical methods 
such that each occurrence of a linear system
can be computed 
through a sequence of
fewer iterations at the expense of not realizing full-precision results, which is then addressed through appropriate levels of preconditioning and iteration for further adjustments.
We summarize such a representative aggressive example in Algorithm~\ref{alg:mkms2} where the RI method is chosen for the iterative approximate numerical method of our general mathematical framework and is solely applied to
\textbf{Step}~$\mathbf{5}$, 
with the understanding that the same approximate mixed-precision numerical method can be 
applied in \textbf{Step}~$\mathbf{3}$;
alternatively, \textbf{Step}~$\mathbf{3}$ can be computed by employing our approximate mixed-precision computing approach in Algorithm~\ref{alg:mkms}.

More formally, define
the block matrices $\bD_{ii} = \bI_{n_i} - \bP_{ii}$, $i\in[m]$;
the strictly block-lower-triangular matrices $\bL_{ij} = \bP_{ij}$, $i > j$, and $\bL_{ij} = \bzero$, $i\leq j$, $i,j\in[m]$;
and the strictly block-upper-triangular matrices $\bUU_{ij} = \bP_{ij}$, $i < j$, and $\bUU_{ij} = \bzero$, $i \geq j$, $i,j \in [m]$.
Further define the block-diagonal matrix $\bD = \diag(\bD_{11},\ldots,\bD_{mm})$. 
We then can decompose the matrix $\bI - \bP$ as
\begin{equation*}
    \bI - \bP = \bD - \bL - \bUU .
\end{equation*}
Now, 
define the diagonal matrix $\bD^{(t-1)}$ of iteration~$t-1$ whose corresponding $i$-th principal diagonal 
block $\bD^{(t-1)}_{ii}$ is set to be $\bD^{(t-1)}_{ii} = (s_i^{(t-1)}/\| \by_i^{(t-1)}\|_1) \, \bI_{n_i}$, $i\in [m]$.
From \textbf{Step}~$\mathbf{4}$, we then can write $\bz^{(t)} = \by^{(t-1)}\bD^{(t-1)}$  which, in combination with the systems of linear equations in \textbf{Step}~$\mathbf{5}$,
leads to $\by^{(t)} = \bz^{(t)}\bL(\bD-\bUU)^{-1}$.
Upon substituting the expression for $\bz^{(t)}$, we finally have
\begin{equation}\label{eq1}
    \by^{(t)} = \by^{(t-1)}\bD^{(t-1)}\bL(\bD-\bUU)^{-1}.
\end{equation}

With this as a starting point, 
the second representative example of a general application of our algorithmic framework within the context
of the KMS algorithm is based on exploiting a combination of~\eqref{eq1} together with levels of mixed-precision preconditioning and iteration for \textbf{Step}~$\mathbf{5}$. 
In particular, transposing both sides of \eqref{eq1} yields
\begin{equation}\label{eq2}
    \by^{(t)\top} = (\bD^{\top}-\bUU^{\top})^{-1} \bL^{\top} 
    \bD^{(t-1)} \by^{(t-1)\top} ,
\end{equation}
from which it is apparent that $\by^{(t)\top}$ can be recast as the solution of a sparse linear system with $\bD^{\top}-\bUU^{\top}$ as its iteration matrix and $\bL^{\top}\bD^{(t-1)} \by^{(t-1)\top}$ as its right-hand side. 
Then the main idea of 
the representative aggressive example of our general algorithmic framework for
\textbf{Step}~$\mathbf{5}$ is to approximately solve the linear systems in~\eqref{eq2} via a fixed number of steps of the mixed-precision 
RI
method~\cite{golub1961chebyshev}. In essence, the preconditioned 
RI
updates the $k_t$-th approximation of $\by^{(t)\top}$ via the fixed-point iteration
\begin{equation} \label{eq3}
    \by^{(t)\top}_{k_t+1}=(\bI-\bM^{-1}(\bD^{\top}-\bUU^{\top}))\by^{(t)\top}_{k_t} + \bM^{-1}\bL^{\top} 
    \bD^{(t-1)} \by^{(t-1)\top},
\end{equation}
where $k_t$ denotes the number of mixed-precision 
RI
steps in the $t$-th iteration of Algorithm~\ref{alg:mkms2} and the matrix $\bM^{-1}$ denotes the preconditioner of the iterative method. For example, if $\bM\equiv \bD^{\top}$, i.e., block-Jacobi preconditioning, the iteration matrix in~\eqref{eq3} becomes equal to $\bD^{-\top}\bUU^{\top}$. Throughout the rest of this paper we consider $\bM$ to be the product of the LU factors of $\bD^{\top}$ obtained using 
reduced precision.\footnote{While not explored in this paper, it is worth mentioning that \eqref{eq3} is equivalent to a linear system with the same iteration matrix $\bD^{\top}-\bUU^{\top}$ and multiple right-hand sides, and thus even further improvements might be possible within the context of our approach; see, e.g.,~\cite{kalantzis2013accelerating,kalantzis2018scalable}.}

\begin{algorithm}[htpb!]
\caption{Representative Aggressive Example of General Algorithmic Framework} \label{alg:mkms2}
\begin{algorithmic}
 \State {\bf Input:} Irreducible stochastic matrix $\bP$ with structural form~\eqref{eq:P-matrix}
  \State {\bf Output:} Estimate of probability distribution $\by$ such that 
  $\by = \by \bP$ and $\by \ones = 1$

  \State Set $t=1$ and Let $\by^{(0)}$ be the given initial approximation to the solution $\by$
  \State Set $\bI-\bP = \bD-\bL-\bUU$ and 
  matrices $\bL_{ij}$ and $\bUU_{ij}$

\State{\bf While:} Stopping criteria is not satisfied
{
\STATE \quad $\mathbf{1.}$ Compute $\hat{\by}_i^{(t-1)} = \frac{\by_i^{(t-1)}}{\| \by_i^{(t-1)} \|_1}$, $i\in[m]$, and Obtain $\hat{\by}^{(t-1)} = (\hat{\by}_1^{(t-1)}, \ldots, \hat{\by}_m^{(t-1)})$
\State \quad $\mathbf{2.}$ Compute $\bQ_{ij}^{(t-1)} = \hat{\by}_i^{(t-1)} \bP_{ij} \ones$, $i,j\in[m]$, and Obtain $\bQ^{(t-1)}$
\State \quad $\mathbf{3.}$ Compute solution of linear system $\bw^{(t-1)} = \bw^{(t-1)} \bQ^{(t-1)}$ and $\bw^{(t-1)} \ones = \| \bw^{(t-1)} \|_1 = 1$ using approximate mixed- \hspace*{0.4in} precision computing method of interest
\State \quad $\mathbf{4.}$ Compute $\bD^{(t-1)}$ whose corresponding $i$-th principal diagonal 
block $\bD^{(t-1)}_{ii}$ is set to be $\bD^{(t-1)}_{ii} = (s_i^{(t-1)}/\| \by_i^{(t-1)}\|_1) \, \bI_{n_i}$, \hspace*{0.4in} $i\in [m]$
\State \quad $\mathbf{5.}$ Compute solution of sparse linear system with $\bD^{\top}-\bUU^{\top}$ as its iteration 
matrix and $\bL^{\top}\bD^{(t-1)} \by^{(t-1)\top}$ as its \hspace*{0.4in} right-hand side using approximate mixed-precision computing method of 
RI in~\eqref{eq3}, 
and Obtain $\by^{(t)}$
\State \quad $\mathbf{6.}$ Conduct convergence test, Setting $t=t+1$ if stopping criteria not satisfied
}

\State{\bf Return:} Estimate $\by^{(t)}$
\end{algorithmic}
\end{algorithm}

\section{Mathematical Analysis}
\label{sec:analysis}
We now turn to derive a mathematical analysis of the 
two representative examples of our
general 
algorithmic
framework of numerical methods from
the previous section.
Our derivation of
a theoretical
error and convergence analysis is presented first for
the representative conservative example of our general mathematical framework
and then for 
the representative aggressive example of our general mathematical framework,
establishing the guaranteed convergence and the rate of convergence for both 
representative examples.
We conclude this section with a performance analysis of certain aspects of the computational improvements provided by our exploitation of multi-precision arithmetic in both
representative examples of our algorithmic framework.
The conditions and properties expressed in both Section~\ref{sec:prelim} and the introduction are assumed to hold throughout our mathematical analysis.

\subsection{Error and Convergence Analysis of 
Representative Conservative Example}
\label{ssec:analysis-KMS}
Our main theoretical result is that 
Algorithm~\ref{alg:mkms}~---~a representative conservative example of our general algorithmic framework~---~is
guaranteed to converge with
an 
approximation error that decreases by a factor of $O(\epsilon)$ at each iteration.
Recall $\by$ to be the exact stationary distribution of the Markov process $\{ X(s) \, ; \, s \in \Ints_+ \}$ and $\by^{(t)}$ to be the estimate of the exact solution $\by$ in iteration~$t$.
Assuming the estimate $\by^{(t-1)}$ of $\by$ has $O(\varphi)$ accuracy, we first present a lemma that shows the estimate $\bw^{(t-1)}$ of the vector $\bw$ in iteration~$t-1$ to also have $O(\varphi)$ accuracy and further shows the estimate $\bQ^{(t-1)}$ of the matrix $\bQ$ in iteration~$t-1$ to have $O(\varphi\epsilon)$ accuracy.
Recalling $\bu_i$ to be the left eigenvector of $\bP_{ii}$ which corresponds to its dominant eigenvalue $\lambda_{i_1} = 1 - \et_i \epsilon + o(\epsilon)$, $i\in[m]$, the lemma also shows that there is an $O(\epsilon)$ norm difference between the exact stationary distribution $\by$ and the corresponding eigenvectors $\bu_i$ of $\bP_{ii}$ multiplied by an $O(1)$ constant.
\begin{lemma}
    Let $\bw = \bw\bQ$ and $\bw^{(t-1)} = \bw^{(t-1)} \bQ^{(t-1)}$ with $\| \bw \|_1 = \| \bw^{(t-1)} \|_1 = 1$, and suppose $\by^{(t-1)}$ is an $O(\varphi)$ approximation to $\by$, $t \in \Nats$, i.e., $\by^{(t-1)}_i = \by_i + O(\varphi)$.
    Then, we have
    \begin{equation}
        \bw^{(t-1)} - \bw = O(\varphi) , \qquad t \in \Nats ,
        \label{eq:w-approx}
    \end{equation}
    and
    the matrix $\bQ^{(t-1)}$ corresponding to $\by^{(t-1)}$ is an $O(\varphi \epsilon)$ approximation of $\bQ$, where
    \begin{equation*}
        \bQ^{(t-1)}_{ij} = \frac{\by^{(t-1)}_i}{\| \by^{(t-1)}_i \|_1} \bP_{ij} \ones , \qquad \forall i, j \in [m].
    \end{equation*}
Further let $\bu_i$ be the left eigenvector of $\bP_{ii}$ corresponding to the exact invariant probability vector $\by_i$, $i\in[m]$.
Then, there exists a constant $\eb_i = O(1)$ such that
\begin{equation}
\| \by_i - \eb_i \bu_i \|_1 = O(\epsilon) .
\label{eq:4.6}
\end{equation}
\label{lem:4.2+4.3}
\end{lemma}

We now establish our main theoretical results that show Algorithm~\ref{alg:mkms} is guaranteed to converge and does so at a rate that decreases the 
approximation error of the estimate $\by^{(t)}$ of the exact solution $\by$ in each outer-loop iteration~$t$ by a factor of $O(\epsilon)$.
Based on our selection of the computational precision with respect to the condition number and the norm of the matrix of the linear system as defined in Section~\ref{ssec:algs-conservative}, these theoretical results also 
imply
that our use of the 
IR
method in the inner loop of Algorithm~\ref{alg:mkms} is guaranteed to converge and does so at a rate that decreases the corresponding relative approximation error in each inner-loop iteration by a factor of $O(\kappa(\bA))$.
\begin{theorem}
Consider Algorithm~\ref{alg:mkms} employing the mixed-precision 
IR
method under the precision-selection process described above in the inner loop \textbf{Step}~$\mathbf{5}$ and \textbf{Step}~$\mathbf{3}$ 
to compute the solution of the corresponding linear system $\bA \bx = \bb$.
Then, 
supposing $\by^{(t-1)}$ is an $O(\varphi)$ approximation to $\by$, $t \in \Nats$, i.e., $\by^{(t-1)} = \by + O(\varphi)$,
Algorithm~\ref{alg:mkms} converges with 
an 
approximation error in the solution $\by^{(t)}$ at each outer-loop iteration~$t$ that decreases by a factor of $O(\epsilon)$, i.e., $\by^{(t)} = \by + O(\varphi \epsilon)$ with $\epsilon < 1$. 
Moreover, the 
IR
method converges linearly with
an 
approximation error that decreases by a factor of $O(\kappa(\bA))$ in each inner-loop iteration.
\label{thm:lem4.1+thm4.1}
\end{theorem}

The results of Lemma~\ref{lem:4.2+4.3} and Theorem~\ref{thm:lem4.1+thm4.1} are based on the assumption that the systems of linear equations in \textbf{Step}~$\mathbf{5}$ and \textbf{Step}~$\mathbf{3}$ are all solved at the exact accuracy of full precision.
We note that the results in Theorem~\ref{thm:lem4.1+thm4.1} on the outer-loop convergence rate associated with a factor of $O(\epsilon)$ and the inner-loop convergence rate associated with a factor of $O(\kappa(\bA))$ hold in a similar manner for the alternative general classes of approximate mixed-precision computing methods of interest.
We also note that the linear convergence rate factor of $O(\kappa(\bA))$ for the inner-loop iterations is typically much faster than the convergence rate factor of $O(\epsilon)$ for the outer-loop iterations of Algorithm~\ref{alg:mkms}.

\subsection{Error and Convergence Analysis of 
Representative Aggressive Example}
\label{ssec:analysis-alternative}
Our main theoretical result is that 
Algorithm~\ref{alg:mkms2}~---~a representative aggressive example of our general algorithmic framework~---~is
guaranteed to converge with
an 
approximation error that decreases with respect to a tradeoff between an increase in the number of iterations~$t$ of the outer loop for fewer numbers of 
RI
steps $k_t$ in \textbf{Step}~$\mathbf{5}$ and a decrease in the computational complexity of the 
RI
method for smaller $k_t$.
To consider these theoretical aspects of 
our second representative example,
we begin by seeking to understand the effects of stopping the mixed-precision 
RI
method after a certain number of steps, say $k\in \Nats$, as formally expressed in the following lemma.
\begin{lemma}\label{lemr1}
Consider computing the solution of the linear system $(\bD^{\top}-\bUU^{\top})\bx=\bb$. 
Then, for any $k>0$, the absolute error satisfies 
\begin{equation*}
    \|\bx-\bx_{[k]}\| \leq \|\bD^{-\top}\| \dfrac{\|\bD^{-\top}\bUU^{\top}\|^{k+1}}{1-\|\bD^{-\top}\bUU^{\top}\|}\|\bb\|.
\end{equation*}
\end{lemma}

We next consider an idealized scenario where $\mathbf{D}^{(t-1)}\equiv \mathbf{C}$, $\forall t\in\Nats$, and $\mathbf{C}$ is a non-singular matrix having the same dimensions as the matrix $\mathbf{D}$. 
The update formula can then be written as $\by^{(t)\top} = ((\mathbf{D}^\top-\mathbf{U}^{\top})^{-1}\mathbf{L}^\top \mathbf{C})
\by^{(t-1)\top}$ from which it follows 
\begin{equation*}
    \by^{(t)\top} = ((\mathbf{D}^\top-\mathbf{U}^{\top})^{-1}\mathbf{L}^\top \mathbf{C})^t\by^{(0)\top} .
\end{equation*}
Assuming that the dominant eigenvalue of the matrix 
$(\mathbf{D}^\top-\mathbf{U}^{\top})^{-1}\mathbf{L}^\top \mathbf{C}$ is simple and has modulus 
$1$,
the iterate $\by^{(t)}$ converges towards the corresponding eigenvector
as formally expressed in the next lemma. 
\begin{lemma} \label{lemr2}
Consider the update $\overline{\by}^{(t)\top} = ((\mathbf{D}^\top-\mathbf{U}^{\top})_{[k_t]}^{-1}\mathbf{L}^\top \mathbf{C})\overline{\by}^{(t-1)\top}$ where $(\mathbf{D}^\top-\mathbf{U}^{\top})^{-1}$ is applied via $k_t$ steps of the 
RI
method, and $k_t$ is monotonically increasing with the number of outer-loop iterations~$t$. 
Then, the sequence $\{\overline{\by}^{(t)}\}$ converges to the same limit as the sequence $\{\by^{(t)}\}$.
\end{lemma}

Our Algorithm~\ref{alg:mkms2} addresses a fundamental tradeoff between computational times on the one hand and error and convergence on the other hand.
At one extreme, if
$k_t$ grows too quickly with the outer-loop iterations $t\in\Nats$, then there will be error and convergence guarantees analogous to Theorem~\ref{thm:lem4.1+thm4.1}, but at the expense of little or no improvement in the computation times.
At the other extreme, if
$k_t$ grows too slowly with the outer-loop iterations $t\in\Nats$, then there will be significant improvements in the computation times, but at the expense of slow convergence and possibly even a lack of guaranteed convergence.
The key to our Algorithm~\ref{alg:mkms2} is the enabling of a control sequence $k_1, k_2, \ldots$ that balances this fundamental tradeoff between the improvements in computation times and the speed of convergence.

The error $\mathbf{e}_k$ after $k$ steps of
RI
satisfies the 
equation $\mathbf{e}_k = (\mathbf{I}-\mathbf{D}^\top+\mathbf{U}^{\top})^k\mathbf{e}_0$, 
where $\mathbf{e}_0$ denotes the initial approximation error. Let us now focus on 
the linear system $\bD^{\top}-\bUU^{\top} \by^{(t)\top}=\bL^{\top}\bD^{(t-1)} 
\by^{(t-1)\top}$. The 
RI method
converges to the true solution 
$\by^{(t)\top}$ as long as the spectral radius
$\rho(\cdot)$
of the matrix $\mathbf{I}-\mathbf{D}^\top+\mathbf{U}^{\top}$ satisfies $\rho(\mathbf{I}-\mathbf{D}^\top+\mathbf{U}^{\top})<1$ 
which is always the case when $\|\mathbf{I}-\mathbf{D}^\top+\mathbf{U}^{\top}\|<1$ 
since $\rho(\mathbf{A}) \leq \|\mathbf{A}\|$ for any matrix $\mathbf{A}$. Therefore, 
after $k_t$ steps of 
RI,
the norm of the approximation error 
$\by^{(t)\top}-\overline{\by}^{(t)\top}$ 
of Lemma~\ref{lemr2}
is reduced by a factor of 
$\rho(\mathbf{I}-\mathbf{D}^\top+\mathbf{U}^{\top})^{k_t}$ compared to the initial error. 
Hence, assuming that we know $\rho(\mathbf{I}-\mathbf{D}^\top+\mathbf{U}^{\top})$, 
we can reduce the initial error by a factor of $\hat{\epsilon}_t> 0$ if we allow 
$k_t$ to be large enough. In particular, assuming that the initial approximation is 
always zero, this leads to 
\begin{equation*}
k_t \geq \dfrac{\mathrm{ln}(\hat{\epsilon}_t/\|\by^{(t)\top}\|)}{\mathrm{ln}(\rho(\mathbf{I}-\mathbf{D}^\top+\mathbf{U}^{\top}))}.
\end{equation*}

The above inequality tells us that, at iteration $t$ of Algorithm~\ref{alg:mkms2}, we can
determine the smallest number of iterations of RI in \textbf{Step}~$\mathbf{5}$ that are required to achieve an arbitrary level of accuracy.
For example,
as used in our numerical experiments,
we assume that $k_1$ is chosen as above and $k_t=2^tk_0$, i.e., the 
number of inner iterations increases geometrically. Starting with an inexact application of the 
RI
method and increasing the number of 
RI
steps $k_t$ with the number of outer-loop iterations $t$,
Algorithm~\ref{alg:mkms2} is able
to avoid over-solving in the earlier stages of \textbf{Step}~$\mathbf{5}$ where the approximate stationary distribution $\by^{(t)}$ is far from $\by$. 
We note here that such tradeoffs are well known in numerical linear algebra and arise often when iterative solvers are used in shift-and-invert algorithms for eigenvalue computations \cite{freitag2007convergence}. 

We now establish our main result for Algorithm~\ref{alg:mkms2} by considering the case where the diagonal matrix $\mathbf{D}^{(t-1)}$ depends (non-linearly) on $\by^{(t-1)}$, which leads to couplings in the dependency of $\by^{(t)}$ with respect to $\by^{(t-1)}$. 
Note that this is exactly the scenario of Algorithm~\ref{alg:mkms2}.
\begin{theorem}\label{thm:mkms2}
Consider Algorithm~\ref{alg:mkms2} employing the mixed-precision 
RI
method under the $k_t$-selection process described above for the number of steps performed in the inner loop \textbf{Step}~$\mathbf{5}$.
Then, supposing $\by^{(t-1)}$ is an $O(\varphi)$ approximation to $\by$, $t \in \Nats$, i.e., $\by^{(t-1)} = \by + O(\varphi)$, 
Algorithm~\ref{alg:mkms2} converges with 
an 
approximation error in the solution $\by^{(t)}$ at each outer-loop iteration~$t$ that decreases by a factor of $O(\epsilon)$, i.e., $\by^{(t)} = \by + O(\varphi \epsilon)$ with $\epsilon < 1$,
and the inner loop of Algorithm~\ref{alg:mkms2} at each outer-loop iteration $t$ converges to within an accuracy of $O(\epsilon^t)$.
Moreover, 
the RI
method converges linearly with
an
approximation error that decreases by a factor of $O(\rho(\mathbf{I}-\mathbf{D}^\top+\mathbf{U}^{\top}))$ in each inner-loop iteration.
\end{theorem}

We note that the results in Theorem~\ref{thm:mkms2} on the outer-loop convergence rate associated with a factor of $O(\epsilon)$ and the inner-loop convergence rate associated with a factor similar to $O(\rho(\mathbf{I}-\mathbf{D}^\top+\mathbf{U}^{\top}))$ can hold in a related manner for appropriate alternatives in the general classes of approximate mixed-precision computing methods of interest.
We also note that the linear convergence rate factor of $O(\rho(\mathbf{I}-\mathbf{D}^\top+\mathbf{U}^{\top}))$ for the inner-loop iterations is typically much faster, or can be made much faster by techniques such as preconditioning, than the convergence rate factor of $O(\epsilon)$ for the outer-loop iterations of Algorithm~\ref{alg:mkms2}.

\subsection{Mixed-Precision Computation Performance Analysis}
\label{sec:performance}
Our general 
mathematical framework of numerical methods exploits
a combination of
advances in mixed-precision computer architectures
and
advances in
iterative approximate computing 
approaches.
Although
well studied in machine learning and related applications, we note that
such
advances in mixed-precision computation have received far less consideration in the context of computing the stationary distribution of Markov processes.
Beyond significantly reducing computation times, our use of lower-precision computation makes it possible to handle much larger systems of linear equations before hitting the memory bandwidth limitations of today's highly-efficient advanced processors.

More specifically, the performance tradeoff at the heart of our 
general
algorithmic
framework
concerns, on the one hand, the significant reductions in execution times afforded by mixed-precision computation at the expense of inaccuracies in the results and, on the other hand, the ability to mitigate inaccurate computations, further reduce execution times and guarantee convergence afforded by iterative approximate computing methods.
Sections~\ref{ssec:analysis-KMS} and~\ref{ssec:analysis-alternative} address the degree to which 
various representative examples of
our algorithmic
framework
mitigate the impact of inaccurate computations, reduce execution times, and ensure convergence.
The reductions in execution times
provided by mixed-precision computation are also quite significant and vary from one computer processor architecture to the next, with increasing trends 
in execution-time reductions from mixed-precision computation
over time.
Generally speaking, there is typically a factor of $2\times$ reduction in computation times with respect to the number of operations per second going from $64$-bit precision to $32$-bit precision, with considerably greater factors of reduction going to $16$-bit precision and to $8$-bit precision including as much as orders of magnitude factors of reduction; 
see, e.g.,~\cite{nvidia1,nvidia2,nvidia3}.
Then, since the cost of any additional iterations incurred by our
representative examples
is often minimal in comparison with the significant reductions from mixed-precision computation 
and iterative approximate approaches
in each iteration as shown in the above theoretical results and the empirical results in Section~\ref{sec:experiments}, our general 
mathematical framework of numerical methods provides
tremendous performance improvements over the most-efficient existing methods.

It is well understood that the computational performance of the advanced processor architectures of today can be hindered by memory bandwidth bottlenecks when the problem is sufficiently large relative to the available memory of the processor.
In particular, rather than realizing the number of operations per second available from the processor in such cases, the computational performance is significantly reduced and instead dictated by the memory bandwidth of the processor architecture.
Therefore, by exploiting reduced-precision computation, our general 
algorithmic framework further enables
us to handle significantly larger systems of linear equations in each iteration at the computational performance afforded by the processor architecture, including the mixed-precision reduction factors above.
This is well beyond what is possible with existing methods whose performance would become memory-bandwidth limited for much smaller problem instances.

\section{Proofs of Main Results}
\label{sec:proofs}
In this section, we provide the proofs of our main theoretical results presented in Section~\ref{sec:analysis}.
The presentations of several of our proofs are 
shortened
with complete proofs of such cases presented in Appendix~\ref{app:proofs}.

\subsection{Proof of Lemma~\ref{lem:4.2+4.3}}
\begin{proof}
We exclude some of the steps and technical details of our proof  
to simplify the presentation,
referring the reader to Appendix~\ref{app:lem:4.2+4.3} for such additional details.

    From the specific supposition $\by^{(t-1)}_i = \by_i + O(\varphi)$ together with $0< \tau \leq \|\by_i\|_1 < 1$, we obtain
    \begin{equation*}
        \| \by^{(t-1)}_i \|_1 = \| \by_i \|_1 + O(\varphi)
        \qquad \mbox{ and } \qquad
        \frac{\| \by^{(t-1)}_i \|_1}{\| \by_i \|_1} = \frac{\| \by_i \|_1}{\| \by_i \|_1} + \frac{O(\varphi)}{\| \by_i \|_1} = 1 + O(\varphi) ,
    \end{equation*}
    from which it follows that
    \begin{equation*}
        \frac{\by^{(t-1)}_i}{\| \by^{(t-1)}_i \|_1} = \frac{\by_i}{\| \by_i \|_1} + O(\varphi) , \qquad i \in [m] .
    \end{equation*}
    This together with the definition of $\bQ^{(t-1)}_{ij} := ({\by^{(t-1)}_i}/{\| \by^{(t-1)}_i \|_1}) \bP_{ij} \ones$ from our Algorithm~\ref{alg:mkms} yields
    \begin{equation*}
        \bQ^{(t-1)}_{ij} = 
        \frac{\by_i}{\| \by_i \|_1} \bP_{ij} \ones + O(\varphi) \bP_{ij} \ones , \qquad \forall i,j \in [m], \; i \neq j.
    \end{equation*}
    Then, from the definition of $\bQ_{ij} := ({\by_i}/{\| \by_i \|_1}) \bP_{ij} \ones$ in the rightmost equation of~\eqref{eq:2.6} and since $\bP_{ij} = O(\epsilon)$ from~\eqref{eq:2.2} for $i \neq j$, we have
    \begin{equation*}
        \bQ^{(t-1)}_{ij} = \bQ_{ij} + O(\varphi \epsilon) , \qquad \forall i,j \in [m], \; i \neq j ,
    \end{equation*}
    which renders the desired result for $\bQ^{(t-1)}$ as an $O(\varphi \epsilon)$ approximation of $\bQ$ when $i \neq j$.
    Given that $\bQ^{(t-1)}$ and $\bQ$ are both stochastic matrices, we further conclude
    \begin{equation*}
        \bQ^{(t-1)}_{ii} = 1 - \sum_{j \in [m] : j \neq i} \bQ^{(t-1)}_{i j} = 1 - \sum_{j \in [m] : j \neq i} \bQ_{i j} + O(\varphi \epsilon) = \bQ_{ii} + O(\varphi \epsilon) , \qquad \forall i \in [m] ,
    \end{equation*}
    thus completing the proof of $\bQ^{(t-1)}$ as an $O(\varphi \epsilon)$ approximation of $\bQ$.

    Next, from the middle equation in~\eqref{eq:4.4}, we have
    \begin{equation}
        (\bI - \bQ) = \bW \begin{pmatrix} 0 & 0 \\ 0 & \bI - \bJ \end{pmatrix} \bW^{-1} .
    \label{eq:4.13}
    \end{equation}
    Given the theorem hypotheses $\bw^{(t-1)} = \bw^{(t-1)} \bQ^{(t-1)}$ and $\bw = \bw \bQ$, 
    where the latter follows from~\eqref{eq:2.8},
    we take the difference between the former
    equation
    and the latter
    equation
    on both sides
    and derive
    \begin{align}
        (\bw^{(t-1)} - \bw) (\bI - \bQ) & = \bw^{(t-1)} ( \bQ^{(t-1)} - \bQ ) \label{eq:4.14} \\
        (\bw^{(t-1)} - \bw) \bW \begin{pmatrix} 0 & 0 \\ 0 & \bI - \bJ \end{pmatrix} & \stackrel{(a)}{=} \bw^{(t-1)} ( \bQ^{(t-1)} - \bQ ) \bW \nonumber,
    \end{align}
    where
    (a)
    follows upon substituting \eqref{eq:4.13} into \eqref{eq:4.14} and multiplying both sides on the right by $\bW$.
    Given the theorem hypothesis $\| \bw \|_1 = \| \bw^{(t-1)} \|_1 = 1$ (according to the rightmost equation in~\eqref{eq:2.8}) and $\| \bW \| = O(1)$ from the rightmost equation in~\eqref{eq:4.4}, together with $\bQ^{(t-1)} = \bQ_{ij} + O(\varphi \epsilon)$ established above, we obtain
    \begin{equation}
        (\bw^{(t-1)} - \bw) \bW \begin{pmatrix} 0 & 0 \\ 0 & \bI - \bJ \end{pmatrix} = \bQ^{(t-1)} - \bQ = O(\varphi \epsilon) .
    \label{eq:4.17}
    \end{equation}
    Letting $(\hat{\ew}_1, \hat{\ew}_2, \ldots, \hat{\ew}_m) = (\bw^{(t-1)} - \bw) \bW$, we have
    \begin{align}
        \hat{\ew}_1 = ((\bw^{(t-1)} - \bw) \bW)_1 = (\bw^{(t-1)} - \bw)\ones & = \| \bw^{(t-1)} \|_1 - \| \bw \|_1 = 0 , \label{eq:4.18} \\
        (\hat{\ew}_2, \ldots, \hat{\ew}_m) (\bI - \bJ) & \stackrel{(a)}{=} O(\varphi \epsilon) , \nonumber
    \end{align}
    where
    (a)
    follows from \eqref{eq:4.17}.
    Since $(\bI - \bJ)$ is nonsingular from~\eqref{eq:4.5}, we first multiply both sides of the last equation on the right by $(\bI - \bJ)^{-1}$ and then take the norm on both sides to conclude
    \begin{equation*}
        \| (\hat{\ew}_2, \ldots, \hat{\ew}_m) \| \leq O(\varphi \epsilon) \| (\bI - \bJ)^{-1} \| .
    \end{equation*}
    From~\eqref{eq:4.5},
    we obtain
    \begin{equation}
        (\hat{\ew}_2, \ldots, \hat{\ew}_m) = O(\varphi) ,
    \label{eq:4.20}
    \end{equation}
    and upon combining \eqref{eq:4.18} and \eqref{eq:4.20}, we have
    \begin{equation*}
        (\bw^{(t-1)} - \bw) \bW = O(\varphi) .
    \end{equation*}
    It therefore follows from the last equation in~\eqref{eq:4.4} that $(\bw^{(t-1)} - \bw) = O(\varphi)$, thus completing the proof of~\eqref{eq:w-approx}.

Lastly, 
from~\eqref{eq:2.4}, 
we have
$\by_i = \sum_{j=1}^m \by_j \bP_{ji}$
and
$\by_i (\bI - \bP_{ii}) = \sum_{j=1:j\neq i}^m \by_j \bP_{ji}$.
Similarly, from~\eqref{eq:4.1} and~\eqref{eq:4.1:cond3}, we obtain
\begin{equation}
    \bu_i \bP_{ii} = (1 - \et_i \epsilon + o(\epsilon) ) \bu_i
    \qquad \mbox{ and } \qquad
     \bu_i ( \bI - \bP_{ii} ) = (\et_i \epsilon + o(\epsilon) ) \bu_i ,
\label{eq:cancel}
\end{equation}
and therefore we conclude
\begin{align}
    \by_i - \bu_i & = \left( \sum_{j=1:j\neq i}^m \by_j \bP_{ji} + ( -\et_i \epsilon + o(\epsilon) ) \bu_i \right) (\bI - \bP_{ii})^{-1} .
    \label{eq:4.7}
\end{align}
Given that
$\begin{pmatrix} \bu_i \\ \bU_i \end{pmatrix}$
is a nonsingular matrix of order $n_i$ and $\sum_{j=1:j\neq i}^m \by_j \bP_{ji}$ is a row vector of order $n_i$, 
there exists an order $n_i$ row vector
\begin{equation}
\bc_i = \sum_{j=1:j\neq i}^m \by_j \bP_{ji}\begin{pmatrix} \bu_i \\ \bU_i \end{pmatrix}^{-1} ,
\label{eq:4.8}
\end{equation}
from which it follows that $\bc_i = O(\epsilon)$ since the norm of $\begin{pmatrix} \bu_i \\ \bU_i \end{pmatrix}^{-1}$ is bounded by a constant independent of $\epsilon$
and $\sum_{j=1:j\neq i}^m \by_j \bP_{ji}$ is $O(\epsilon)$.
Then,
upon substitution of~\eqref{eq:4.8} into~\eqref{eq:4.7}, we obtain
\begin{align}
\by_i - \bu_i 
& \stackrel{(a)}{=} \bc_i \begin{pmatrix} (\et_i\epsilon + o(\epsilon))^{-1} & 0 \\ \bzero & (\bI - \bH_i)^{-1} \end{pmatrix} \begin{pmatrix} \bu_i \\ \bU_i \end{pmatrix} - \bu_i \nonumber \\
&= \frac{\ec_{i_1}}{\et_i\epsilon + o(\epsilon)} \bu_i + (\ec_{i_2}, \ldots, \ec_{i_{n_i}}) (\bI - \bH_i)^{-1} \bU_i - \bu_i ,
\label{eq:4.8+}
\end{align}
where (a) follows upon substituting~\eqref{eq:4.1:cond3} together with~\eqref{eq:cancel}. 
Given that $\bc_i = O(\epsilon)$, $(\bI - \bH_i)^{-1}$ is $O(1)$ from~\eqref{eq:4.2}
and $\bU_i$ is also $O(1)$ from above, it follows from~\eqref{eq:4.8+} that
\begin{align*}
    \by_i - \frac{\ec_{i_1}}{\et_i\epsilon + o(\epsilon)} \bu_i & = O(\epsilon) .
\end{align*}
Let $\eb_i = \ec_{i_1}/(\et_i\epsilon + o(\epsilon))$.
Since $\bc_i = O(\epsilon)$, this implies that $\ec_{i_1}$ is of order $\es_i \epsilon + o(\epsilon)$ for some $\es_i > 0$, and thus we have
\begin{equation*}
    \eb_i = \frac{\ec_{i_1}}{\et_i\epsilon + o(\epsilon)} = \frac{\es_i \epsilon + o(\epsilon)}{\et_i\epsilon + o(\epsilon)} = \frac{\es_i}{\et_i} + o(1) = O(1) ,
\end{equation*}
from which the desired result~\eqref{eq:4.6} follows.
\end{proof}

\subsection{Proof of Theorem~\ref{thm:lem4.1+thm4.1}}
\begin{proof}
We exclude some of the steps and technical details of our proof 
to simplify the presentation,
referring the reader to Appendix~\ref{app:thm:lem4.1+thm4.1} for such additional details.

Let
the error in the approximate solution $\by^{(t-1)}$ at iteration~$t-1$ of Algorithm~\ref{alg:mkms} be of order $\varphi$, namely $\by^{(t-1)} = \by + O(\varphi)$.
Then, from Lemma~\ref{lem:4.2+4.3}, we have
\begin{equation}
    \bQ^{(t-1)} = \bQ + O(\varphi \epsilon) \qquad \mbox{and} \qquad
    \bw^{(t-1)} = \bw + O(\varphi) .
\label{eq:lemma4.1}
\end{equation}
Starting with the $m$-th component of $\by^{(t)}$ for iteration~$t$ of Algorithm~\ref{alg:mkms}, we derive
\begin{align}
    \by_m^{(t)} & = \by_m^{(t)} \bP_{mm} + \sum_{j=1}^{m-1} s_j^{(t-1)} \hat{\by}_j^{(t-1)} \bP_{jm} \nonumber \\
    \by_m^{(t)} (\bI - \bP_{mm}) & \stackrel{(a)}{=} \sum_{j=1}^{m-1} (s_j + O(\varphi)) (\hat{\by}_j + O(\varphi)) \bP_{jm} \nonumber \\
    & \stackrel{(b)}{=} \by_m (\bI - \bP_{mm}) + \sum_{j=1}^{m-1} O(\varphi) \bP_{jm} ,
    \label{eq:4.25}
\end{align}
where
(a) follows from substituting the rightmost equation in~\eqref{eq:lemma4.1} and from the $O(\varphi)$ approximation of $\by^{(t-1)}$,
and
(b) follows from~\eqref{eq:2.6} and~\eqref{eq:2.9} and from $s_j$ and $\hat{\by}_j$ both being $O(1)$ together with $\varphi < 1$.
Similarly to~\eqref{eq:4.8},
given that
$\begin{pmatrix} \bu_m \\ \bU_m \end{pmatrix}$
is a nonsingular matrix of order $n_m$
whose norm is bounded by a constant independent of $\epsilon$,
and given that $\sum_{j=1}^{m-1} O(\varphi) \bP_{jm}$ is a row vector of order $n_m$ whose components are $O(\varphi\epsilon)$, there exists a row vector $\bar{\bc}_m$ of order $n_m$ whose components are $O(\varphi\epsilon)$ such that
\begin{equation}
\sum_{j=1}^{m-1} O(\varphi) \bP_{jm} = \bar{\bc}_m \begin{pmatrix} \bu_m \\ \bU_m \end{pmatrix} .
\label{eq:4.26}
\end{equation}
Multiplying both sides of~\eqref{eq:4.25} by $(\bI - \bP_{mm})^{-1}$ and substituting~\eqref{eq:4.26}, we obtain
\begin{align}
    \by_m^{(t)} - \by_m 
    & \stackrel{(a)}{=} \bar{\bc}_m
    \begin{pmatrix} \displaystyle \frac{\bu_m}{\et_m\epsilon + o(\epsilon)} \vspace*{0.15cm} \\ (\bI - \bH_m)^{-1} \bU_m \end{pmatrix} \nonumber \\
    & = \left( \frac{\bar{\ec}_{m_1}}{\et_m\epsilon + o(\epsilon)} \right) \bu_m + (\bar{\ec}_{m_2}, \ldots, \bar{\ec}_{m_{n_m}}) (\bI - \bH_m)^{-1} \bU_m ,
    \label{eq:4.27}
\end{align}
where (a) follows upon substitution of~\eqref{eq:4.1:cond3}.
From~\eqref{eq:4.6}, we have
\begin{equation}
    \bu_m = \eb_m^{-1} \by_m + O(\epsilon) ,
    \label{eq:4.28}
\end{equation}
which together with~\eqref{eq:4.27}, \eqref{eq:4.28}, $\bar{\ec}_{m_k} = O(\varphi\epsilon)$, $\ec_{m_1} = O(\epsilon)$, $(\bI - \bH_m)^{-1} = O(1)$ and $\bU_m = O(1)$ yields
\begin{align*}
    \by_m^{(t)} - \by_m & = \left( \frac{\bar{\ec}_{m_1}}{\et_m\epsilon + o(\epsilon)} \right) \left(\frac{\et_m\epsilon + o(\epsilon)}{\ec_{m_1}} \by_m + O(\epsilon) \right) + (\bar{\ec}_{m_2}, \ldots, \bar{\ec}_{m_{n_m}}) (\bI - \bH_m)^{-1} \bU_m \\
    \by_m^{(t)} & = (1 + O(\varphi)) \by_m + O(\varphi \epsilon) .
\end{align*}
It then follows that
\begin{equation*}
    \hat{\by}_m^{(t)} = \frac{\by_m^{(t)}}{\| \by_m^{(t)} \|_1} = \frac{(1 + O(\varphi)) \by_m + O(\varphi \epsilon)}{\| (1 + O(\varphi)) \by_m + O(\varphi \epsilon) \|_1} = \frac{\by_m}{\| \by_m \|_1} + O(\varphi \epsilon) ,
\end{equation*}
from which we conclude
\begin{equation}
    \hat{\by}_m^{(t)} = \hat{\by}_m + O(\varphi \epsilon) .
    \label{eq:4.30}
\end{equation}

Now, arguing by induction with respect to the base case~\eqref{eq:4.30}, 
we can therefore establish
\begin{align*}
    \frac{\by_i^{(t)}}{\| \by_i^{(t)} \|_1} & = \frac{\by_i}{\| \by_i \|_1}+ O(\varphi \epsilon) , \qquad \forall i \in [m] .
\end{align*}
It then follows from the initial arguments in the proof of Lemma~\ref{lem:4.2+4.3} that
\begin{equation*}
    \by^{(t)}_i = \by_i + O(\varphi \epsilon) , \qquad \forall i \in [m] ,
\end{equation*}
which yields the desired result.

Next, 
under the process employed for selecting the precision of computations in \textbf{Step}~$\mathbf{5}$ and \textbf{Step}~$\mathbf{3}$ of Algorithm~\ref{alg:mkms} with respect to the condition number and the norm of the matrix $\bA$ of the linear system, we know that the assumptions of the mixed-precision analysis of 
the IR method
by Moler~\cite{Moler} are satisfied.
It then follows from the convergence results of Moler~\cite{Moler} that, under these conditions, the
IR
method is guaranteed to converge linearly with an approximation error that decreases by a factor of $O(\kappa(\bA))$ in each iteration, thus completing the proof.
\end{proof}

\subsection{Proof of Lemma~\ref{lemr1}}
\begin{proof}
First, notice that the matrix $\bD^{\top}-\bUU^{\top}$ is block-lower triangular, with the $i$-th diagonal block being equal to $\bI -\bP_{ii}$. Since $\bP_{ii}$ has only one eigenvalue that is close to one, and since this eigenvalue $\lambda_{i_1}$ is equal to $1-\et_i\epsilon + o(\epsilon)$, the smallest eigenvalue of $\bD^{\top}-\bUU^{\top}$ is of the order $O(\epsilon)$.  
Expressing
\begin{align*}
    \bD^{\top}-\bUU^{\top} = \bD^{\top}(\bI-\bD^{-\top}\bUU^{\top}),
\end{align*}
and inverting both sides yields 
\begin{align*}
    (\bD^{\top}-\bUU^{\top})^{-1} = (\bI-\bD^{-\top}\bUU^{\top})^{-1}\bD^{-\top}.
\end{align*}
Assuming the invertibility of $\bD^{\top}-\bUU^{\top}$, it follows that $\|\bD^{-\top}\bUU^{\top}\| <1$, and thus $(\bI-\bD^{-\top}\bUU^{\top})^{-1}$ possesses a convergent Neumann series  that can be expanded as 
\begin{align*}
    (\bI - \bD^{-\top}\bUU^{\top})^{-1} = \bI + \bD^{-\top}\bUU^{\top} + (\bD^{-\top}\bUU^{\top})^2 + (\bD^{-\top}\bUU^{\top})^3 + \cdots,
\end{align*}
which leads to
\begin{align*}
    (\bD^{\top}-\bUU^{\top})^{-1} = \left(\bI + \bD^{-\top}\bUU^{\top} + (\bD^{-\top}\bUU^{\top})^2 + (\bD^{-\top}\bUU^{\top})^3 + \cdots\right)\bD^{-\top}.
\end{align*}

Now, upon truncating the above expression for $(\bD^{\top}-\bUU^{\top})^{-1}$ up to $k$ 
terms, we obtain 
\begin{align*}
    (\bD^{\top}-\bUU^{\top})_{[k]}^{-1} = \left(\bI + \bD^{-\top}\bUU^{\top} + (\bD^{-\top}\bUU^{\top})^2 + (\bD^{-\top}\bUU^{\top})^3 + \cdots + (\bD^{-\top}\bUU^{\top})^k\right)\bD^{-\top},
\end{align*}
which results in an error of
\begin{align*}
    (\bD^{\top}-\bUU^{\top})^{-1} - (\bD^{\top}-\bUU^{\top})_{[k]}^{-1} = \left((\bD^{-\top}\bUU^{\top})^{k+1} + (\bD^{-\top}\bUU^{\top})^{k+2} + \cdots\right)\bD^{-\top}.
\end{align*}
Recalling from above that $\|\bD^{-\top}\bUU^{\top}\| <1$, we obtain the inequality 
\begin{equation*}
    \|(\bD^{\top}-\bUU^{\top})^{-1} - (\bD^{\top}-\bUU^{\top})_{[k]}^{-1}\| \leq \|\bD^{-\top}\| \dfrac{\|\bD^{-\top}\bUU^{\top}\|^{k+1}}{1-\|\bD^{-\top}\bUU^{\top}\|} .
\end{equation*}
From this upper bound on the error we deduce that, when $\|\bD^{-\top}\bUU^{\top}\| \ll 1$, a small value of $k$ can lead to an accurate approximation of the inverse 
$(\bD^{\top}-\bUU^{\top})^{-1}$, which yields the desired result. 
\end{proof}

\subsection{Proof of Lemma~\ref{lemr2}}
\begin{proof}
The expression for the update of the lemma can be equivalently written as 
\begin{equation*}
\overline{\by}^{(t)\top} = \prod_{t=1} (((\bD^{\top}-\bUU^{\top})^{-1} - ((\bD^{-\top}\bUU^{\top})^{k_t+1} + (\bD^{-\top}\bUU^{\top})^{k_t+2} + \cdots)\bD^{-\top}) \mathbf{L}^\top \mathbf{C})\by^{(0)\top}.
\end{equation*}
Then, in order for the sequence $\{\overline{\by}^{(t)}\}$ to converge to the same limit as the sequence 
$\{\by^{(t)}\}$, we need to have
$$\prod_{t=1} (((\bD^{\top}-\bUU^{\top})^{-1} - ((\bD^{-\top}\bUU^{\top})^{k_t+1} + (\bD^{-\top}\bUU^{\top})^{k_t+2} + \cdots)\bD^{-\top})\mathbf{L}^\top \mathbf{C})$$ 
converge to the dominant eigendirection of $(\mathbf{D}^\top-\mathbf{U}^{\top})^{-1}\mathbf{L}^\top \mathbf{C}$, which is true provided that $k_t$ increases monotonically with the number of outer-loop iterations~$t$, thus completing the proof.
\end{proof}

\subsection{Proof of Theorem~\ref{thm:mkms2}}
\begin{proof}
Let
$F$ denote the (fixed) process defined by a single outer iteration of Algorithm~\ref{alg:mkms2},
and thus ${\by}^{(t)\top} = F({\by}^{(t-1)\top})$. 
By supposition, we have 
\begin{equation}\label{eq:thm2:supposition}
    \|{\by}^{(t-1)\top} -{\by}^{\top}\| = O(\varphi) .
\end{equation}
Now, instead of ${\by}^{(t)\top}$, assume that Algorithm~\ref{alg:mkms2} computes  
\begin{equation}\label{eq:forward-error}
    {\by}^{(t)\top} = F({\by}^{(t-1)\top}) + \mathbf{r}_t^{k_t}, 
\end{equation}
where $\mathbf{r}_t^{k_t}$ denotes the approximation error resulting from the application of 
$k_t$ steps of
the RI method.
Subtracting ${\by}^{\top}$ from both sides of~\eqref{eq:forward-error},
we have 
\begin{equation*}
    {\by}^{(t)\top} -{\by}^{\top} = F({\by}^{(t-1)\top})-{\by}^{\top} + \mathbf{r}_t^{k_t},
\end{equation*}
and thus taking the norms of both sides and substituting for $\mathbf{r}_t^{k_t}$ yields
\begin{equation*}
    \|{\by}^{(t)\top} -{\by}^{\top}\| \leq \|F({\by}^{(t-1)\top})-{\by}^{\top}\| + \rho(\mathbf{I}-\mathbf{D}^\top+\mathbf{U}^{\top})^{k_t} \|\mathbf{L}^\top \mathbf{D}^{(t-1)}{\by}^{(t-1)\top}\|.
\end{equation*}
It then follows from analogous arguments establishing Theorem~\ref{thm:lem4.1+thm4.1} that
$\|F({\by}^{(t-1)\top})-{\by}^{\top}\| = O(\epsilon)\|{\by}^{(t-1)\top} -{\by}^{\top}\|$, and therefore we can simplify the above equation as
\begin{equation*}
    \|{\by}^{(t)\top} -{\by}^{\top}\| \leq O(\epsilon)\|{\by}^{(t-1)\top} -{\by}^{\top}\| + \rho(\mathbf{I}-\mathbf{D}^\top+\mathbf{U}^{\top})^{k_t} \|\mathbf{L}^\top \mathbf{D}^{(t-1)}{\by}^{(t-1)\top}\|.
\end{equation*}
Recall though that $k_t$ is
always set
so that 
$\|\mathbf{r}_t^{k_t}\|\leq 
O(\varphi\epsilon)$.
Therefore, we can replace 
$\rho(\mathbf{I}-\mathbf{D}^\top+\mathbf{U}^{\top})^{k_t} \|\mathbf{L}^\top \mathbf{D}^{(t-1)}{\by}^{(t-1)\top}\|$ by its upper-bound, leading to 
\begin{align*}
    \|{\by}^{(t)\top} -{\by}^{\top}\| &\leq O(\epsilon)\|{\by}^{(t-1)\top} -{\by}^{\top}\| + 
    O(\varphi\epsilon)\\
    &\leq O(\varphi\epsilon) ,
\end{align*}
where the last result follows from~\eqref{eq:thm2:supposition} and $\epsilon<1$.
Finally, the fact that RI converges linearly with respect to 
$\rho(\mathbf{I}-\mathbf{D}^\top+\mathbf{U}^{\top})$ follows by recalling that 
$\mathbf{r}_t^{k_t} = (\mathbf{I}-\mathbf{D}^\top+\mathbf{U}^{\top})^{k_t}\mathbf{L}^\top \mathbf{D}^{(t-1)}{\by}^{(t-1)\top}$.
\end{proof}

\section{Numerical Experiments}
\label{sec:experiments}
In this section we present a representative sample of numerical experiments that support our theoretical results and empirically evaluate our general 
mathematical framework of numerical methods
in comparison with the most-efficient existing 
algorithms
for 
computing the stationary distribution of
NCD Markov processes.
The
most-efficient existing algorithms are the iterative numerical methods due to Takahashi~\cite{Taka75}, Vantilborgh~\cite{Vant81}, and Koury, McAllister and Stewart~\cite{KoMcSt84}
(KMS),
all three 
of which
exploit aggregation-disaggregation 
in a very similar manner with similar convergence behaviors~\cite{CaoSte85}.
Although our general algorithmic framework of computational approaches can be applied within the context of any of these existing numerical methods and beyond, 
we again select
the 
KMS
algorithm 
as the appropriate baseline method for the numerical evaluation comparison of 
the two representative examples of our algorithmic framework in Section~\ref{sec:algs}, since KMS is the most recent of the most-efficient existing algorithms.

We first present details on the set up of our numerical experiments and then we present and discuss the performance results from a representative sample of our numerical experiments.
Beyond quantitatively supporting our theoretical results, these empirical performance results demonstrate that 
the two representative examples of
our general 
algorithmic 
framework
exhibit relatively little or no increase in the number of outer-loop iterations and orders of magnitude improvements in the computation time, both in comparison with the baseline method of the
KMS algorithm.

\subsection{Experimental Setup}
Our numerical experiments of 
computing the stationary distribution of
NCD Markov processes with  Algorithm~\ref{alg:mkms}, Algorithm~\ref{alg:mkms2} and the baseline
KMS algorithm
are conducted in Python on a standard compute infrastructure with a $2.5$Ghz i$7$ processor and $100$GB of RAM.
For one primary set of numerical experiments, the
$n\times n$ transition probability matrix $\bP$ of block structural form~\eqref{eq:P-matrix} for each independent trial run is randomly generated with each element sampled independently and identically from a uniform distribution on the unit interval. 
Then, to satisfy the conditions~\eqref{eq:2.2}, we reduce the magnitude of the elements of the off-diagonal block matrices by proportionally scaling these elements such that $\|\bP_{ij}\| = O(\epsilon)$, $i,j\in[m]$, $i\neq j$;
and, to satisfy the requirement that $\bP$ be a stochastic matrix, we correspondingly and proportionally scale each row of the matrices $\bP_{ii}$, $i\in[m]$, so that each row of the resulting matrix $\bP$ sum to one.
This first collection of numerical experiments essentially supports our empirical performance evaluation of Algorithm~\ref{alg:mkms} and Algorithm~\ref{alg:mkms2} with respect to the various parameters of the NCD Markov processes, such as $m$, $\epsilon$ and $n_i$, $i\in[m]$.
For another primary set of numerical experiments, we take a similar approach while exploiting a standard library of matrices for real-world linear systems~\cite{SpareMatrixLibrary}.
In particular, as described above for each independent trial run, the $n\times n$ transition probability matrix $\bP$ is randomly generated with each element sampled independently and identically from a uniform distribution on the unit interval.
Then, the $n_i\times n_i$ diagonal block matrices $\bP_{ii}$, $i\in[m]$, are replaced by real-world matrices selected from the standard library in~\cite{SpareMatrixLibrary}.
Next, to satisfy the conditions~\eqref{eq:2.2}, we reduce the magnitude of the elements of the off-diagonal block matrices by proportionally scaling these elements such that $\|\bP_{ij}\| = O(\epsilon)$, $i,j\in[m]$, $i\neq j$;
and, to satisfy the requirement that $\bP$ be a stochastic matrix, we correspondingly and proportionally scale each row of the matrices $\bP_{ii}$, $i\in[m]$, so that each row of the resulting matrix $\bP$ sum to one.
This second collection of numerical experiments based on real-world matrices essentially supports and validates our first collection of numerical experiments based on randomly generated matrices, both with respect to the empirical performance evaluation of our Algorithm~\ref{alg:mkms} and Algorithm~\ref{alg:mkms2}.

The performance metrics of interest taken over the independent trial runs of our numerical experiments are the mean number of outer-loop iterations and the mean computation time.
To support some of the most advanced GPUs currently available in the marketplace, we developed a simulation of the computation times on a given GPU based on the specifications of the corresponding GPU architecture;
and we validated our simulation of the computation times against physical experiments on the NVIDIA A100 processor~\cite{NvidiaA100}.
More specifically, the computation times for every step of each of the three algorithms are simulated with respect to the number of operations for each step of the corresponding iterative process at its level of precision and the processor performance specification of the GPU of interest.
For every algorithm, 
the computation times
required to 
solve
a $k\times k$ linear system 
are based on the use of
LU decomposition 
at the precision level employed by the given algorithm for this linear system computation.
Other than \textbf{Step}~$\mathbf{5}$ and \textbf{Step}~$\mathbf{3}$, which are performed in lower precision within Algorithm~\ref{alg:mkms} and Algorithm~\ref{alg:mkms2} (see below),
the simulated computation times for each outer-loop iteration are exactly the same for
all three algorithms 
and
are performed in full-precision.
For 
the purpose of our numerical experiments,
we used the performance of the NVIDIA H100 processor~\cite{nvidia1,nvidia2} whose peak performance for double-, single- and half-precision floating point arithmetic are 34, 67 and 134 Tflops per second, respectively.
The baseline 
KMS algorithm performs all of its computations in double-precision;
Algorithm~\ref{alg:mkms} performs the computations of \textbf{Step}~$\mathbf{5}$ and \textbf{Step}~$\mathbf{3}$ based on its precision-selection process, though typically in single-precision;
Algorithm~\ref{alg:mkms2} performs the computations of \textbf{Step}~$\mathbf{5}$ in single-precision and the computations of \textbf{Step}~$\mathbf{3}$ in full-precision.

The 
simulation of
computation times 
described above
are only used when the matrix of each linear system fits within the available memory on the 
GPU
processor architecture.
In particular, when the size of the matrix $\bP_{ii}$, $i\in[m]$, at the level of precision employed exceeds the available memory of $96$GB, then the computation of the linear system becomes memory bandwidth limited as described in Section~\ref{sec:performance}.
Based on the specification for the NVIDIA H100 processor~\cite{nvidia1,nvidia2}, the memory bandwidth of $3.35$ Tbytes per second is instead used when any of the algorithms exceed the available memory and become bandwidth limited.
As noted in Section~\ref{sec:performance}, the use of lower-precision computation in our general algorithmic approaches makes it possible to handle larger systems of linear equations before hitting the memory-bandwidth limitations of the processor architecture.

Lastly, it is important to note that the performance benefits of Algorithm~\ref{alg:mkms} and Algorithm~\ref{alg:mkms2} over the baseline 
KMS algorithm
in sharply reducing the core computational bottlenecks associated with solving systems of linear equations
stem from a combination of the higher performance throughput of the lower-precision arithmetic and the reduced computations of the iterative approximate numerical methods.

\subsection{Empirical Performance Results}
We now present the performance results from our numerical experiments of 
computing the stationary distribution of
NCD Markov processes with Algorithm~\ref{alg:mkms} and Algorithm~\ref{alg:mkms2} that support our theoretical results and 
quantify the empirical benefits of our general 
mathematical framework of numerical methods
over the baseline
KMS algorithm.
The primary parameters that influence the computation of the solution of~\eqref{eq:2.4} to obtain the stationary distribution $\by$ of NCD Markov processes concern:
($a$) the dimensions $n_i$ of the diagonal block matrices $\bP_{ii}$, $i\in[m]$;
($b$) the number $m$ of diagonal block matrices; and
($c$) the magnitude $\epsilon$ of transitions in the off-diagonal block matrices $\|\bP_{ij}\| = O(\epsilon)$, $i,j\in[m], i\neq j$.
We 
first
consider collections of performance result comparisons from numerical experiments in which two of the primary parameters are fixed and the remaining primary parameter varies over a range of values.
Each set of performance result comparisons comprises evaluation of  
the representative conservative example of our general algorithmic framework in Algorithm~\ref{alg:mkms},
the representative aggressive example of our general algorithmic framework in Algorithm~\ref{alg:mkms2},
and the baseline KMS algorithm.
As previously noted, the
quantitative metrics of interest in each collection of performance comparisons are the mean number of iterations and the mean computation time in seconds, taken over ten independent experimental trial runs.

For Algorithm~\ref{alg:mkms}, we select the computational precision level with respect to the condition number and the norm of the matrix of the linear systems of \textbf{Step}~$\mathbf{5}$ and \textbf{Step}~$\mathbf{3}$ as defined in Section~\ref{ssec:algs-conservative}.
We then employ the 
IR
method as described in Section~\ref{ssec:algs-conservative} to compute the solution of each linear system to full-precision accuracy.
For Algorithm~\ref{alg:mkms2}, we employ the 
RI
method with block-Jacobi preconditioning as described in Section~\ref{ssec:algs-aggressive} 
to solve the linear systems of \textbf{Step}~$\mathbf{5}$ with $\bD^{\top}-\bUU^{\top}$ as the iteration matrix and $\bL^{\top}\bD^{(t-1)} \by^{(t-1)\top}$ as the right-hand side.
In particular, the linear systems of \textbf{Step}~$\mathbf{5}$ with the coefficient matrix $\bD^{\top}-\bUU^{\top}$ are preconditioned by the block-diagonal matrix $\bD^{\top}$, where the preconditioner is applied via LU factorization computed in $32$-bit precision. 
For the solution of each linear system of \textbf{Step}~$\mathbf{5}$, we start with a number of maximum-allowed iterations equal to ten and
then geometrically increase this value per each outer iteration. 
The linear system of \textbf{Step}~$\mathbf{3}$ is solved exactly in full precision so that we can focus on the performance benefits of our general algorithmic
framework
for the larger linear systems of \textbf{Step}~$\mathbf{5}$.

\begin{figure}[htbp]
\centering
\begin{subfigure}{0.37\textwidth}
    \includegraphics[width=\textwidth]{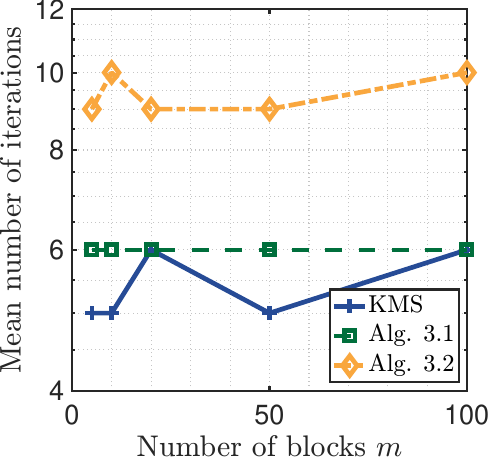}
    \caption{Mean number of iterations.}
    \label{fig1:a}
\end{subfigure}
\hfill
\begin{subfigure}{0.37\textwidth}
    \includegraphics[width=\textwidth]{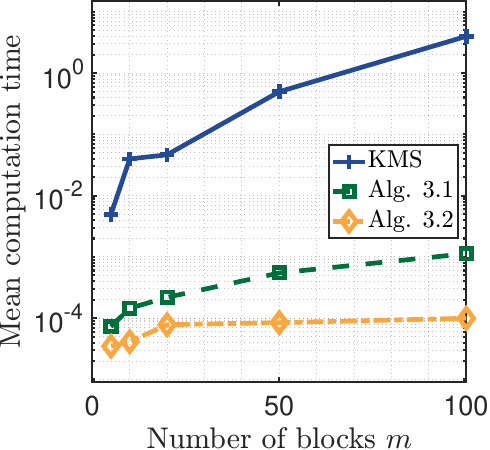}
    \caption{Mean computation times (seconds).}
    \label{fig1:b}
\end{subfigure}
\caption{Performance comparison of our 
Algorithm~\ref{alg:mkms} and Algorithm~\ref{alg:mkms2}
with the KMS baseline for diagonal matrices $\bP_{ii}$ having dimension $n_i=500$ and off-diagonal matrices $\|\bP_{ij}\| = O(\epsilon)$ having $\epsilon = 0.1$, as a function of the number of blocks $m\in\{5, 10, 20, 50, 100\}$.
The comparative performance results represent averages over $10$ independent trial runs executed on a simulated H$100$ GPU architecture.}
\label{fig:1}
\end{figure}

\begin{figure}[htbp]
\centering
\begin{subfigure}{0.37\textwidth}
    \includegraphics[width=\textwidth]{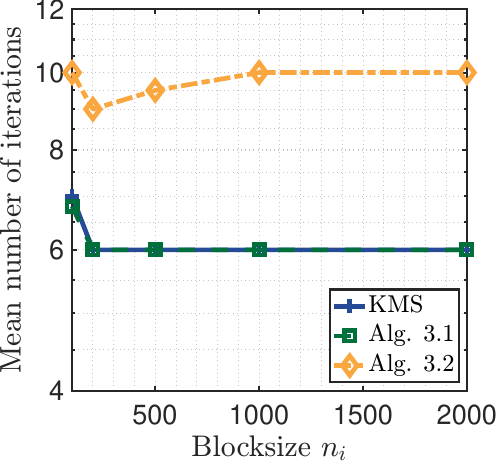}
    \caption{Mean number of iterations.}
    \label{fig2:a}
\end{subfigure}
\hfill
\begin{subfigure}{0.37\textwidth}
    \includegraphics[width=\textwidth]{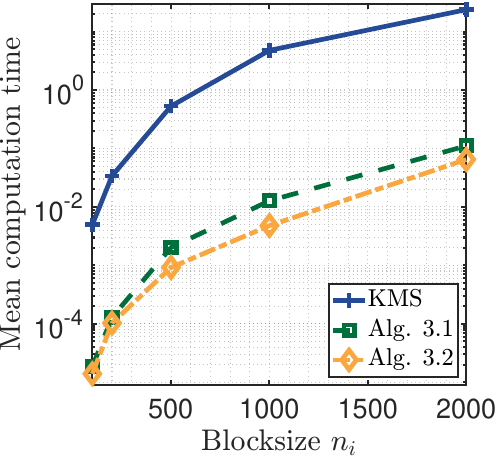}
    \caption{Mean computation times (seconds).}
    \label{fig2:b}
\end{subfigure}
\caption{Performance comparison of our Algorithm~\ref{alg:mkms} and Algorithm~\ref{alg:mkms2}
with the KMS baseline for number of blocks $m=20$ and off-diagonal matrices $\|\bP_{ij}\| = O(\epsilon)$ having $\epsilon = 0.1$, as a function of the dimension
$n_i\in\{100, 200, 500, 1000, 2000\}$
of the diagonal matrices $\bP_{ii}$.
The comparative performance results represent averages over $10$ independent trial runs executed on a simulated H$100$ GPU architecture.}
\label{fig:2}
\end{figure}

\begin{figure}[htbp]
\centering
\begin{subfigure}{0.37\textwidth}
    \includegraphics[width=\textwidth]{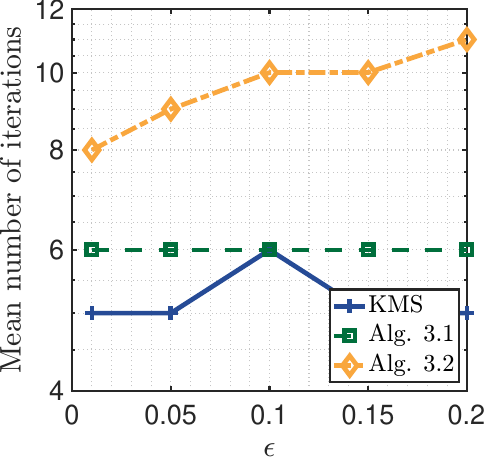}
    \caption{Mean number of iterations.}
    \label{fig3:a}
\end{subfigure}
\hfill
\begin{subfigure}{0.37\textwidth}
    \includegraphics[width=\textwidth]{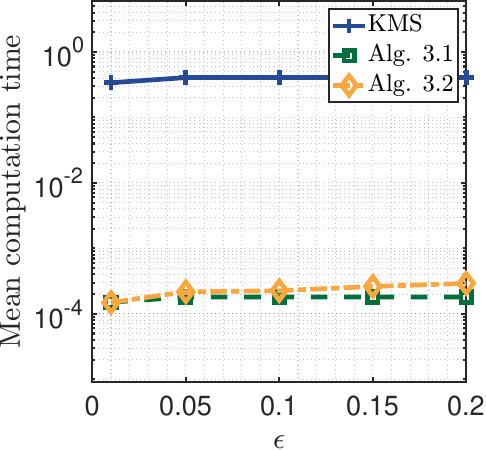}
    \caption{Mean computation times (seconds).}
    \label{fig3:b}
\end{subfigure}
\caption{Performance comparison of our Algorithm~\ref{alg:mkms} and Algorithm~\ref{alg:mkms2}
with the KMS baseline for diagonal matrices $\bP_{ii}$ having dimension $n_i=500$ and number of blocks $m=20$, as a function of $\epsilon\in\{0.01, 0.05, 0.1, 0.15, 0.2\}$ for the off-diagonal matrices $\|\bP_{ij}\| = O(\epsilon)$.
The comparative performance results represent averages over $10$ independent trial runs executed on a simulated H$100$ GPU architecture.}
\label{fig:3}
\end{figure}

We present
in Figures~\ref{fig:1}~--~\ref{fig:3}
the results from a representative sample of our 
first collection of
numerical experiments
based on randomly generated $n\times n$ transition probability matrices $\bP$ under different NCD process parameters.
Figure~\ref{fig:1} presents the representative set of performance result comparisons for the case where the dimensions $n_i$ of the diagonal block matrices $\bP_{ii}$, $i\in[m]$, are fixed to be $500$, the magnitude $\epsilon$ of transitions in the off-diagonal block matrices $\|\bP_{ij}\| = O(\epsilon)$, $i,j\in[m], i\neq j$, is fixed to be $0.1$, and the number $m$ of diagonal block matrices is varied such that $m\in\{5, 10, 20, 50, 100\}$.
We first observe from Figure~\ref{fig:1}(a) that the mean number of iterations under Algorithm~\ref{alg:mkms} are either identical to that of the baseline KMS algorithm or differ by at most one, with the mean number of iterations under Algorithm~\ref{alg:mkms2} only somewhat higher requiring a few additional iterations (though of lower computational complexity).
The mean number of iterations as a function of the number $m$ of block matrices $\bP_{ii}$, $i\in[m]$, remains fairly flat for all algorithms, which is as expected since the linear increase is in the number of linear systems of the same size being solved.
From Figure~\ref{fig:1}(b), we observe that both of 
the representative examples of
our general algorithmic
framework
provide orders of magnitude reduction in the mean computation time and that such performance improvements increase with the number $m$ of block matrices $\bP_{ii}$, $i\in[m]$.
The mean computation time for each algorithm increases with $m$ as expected due to the linear increase in the number of linear systems of the same size being solved, with the rate of increase greatest for the baseline KMS algorithm, with much smaller rates of increase for both of 
the representative examples of
our
algorithmic 
framework,
and with even smaller rates of increase for Algorithm~\ref{alg:mkms2} than for Algorithm~\ref{alg:mkms}.
We note that our choice for the number of maximum-allowed iterations in the solution of each linear system of \textbf{Step}~$\mathbf{5}$ under Algorithm~\ref{alg:mkms2} leads to the lowest mean computation times in Figure~\ref{fig:1}(b).
The increase in the number $m$ of block matrices $\bP_{ii}$, $i\in[m]$, under Algorithm~\ref{alg:mkms2} increases the part of the iteration matrix not captured by the block-Jacobi preconditioner, but this use of the method becomes computationally less expensive and thus there is an interesting algorithmic tradeoff 
in Algorithm~\ref{alg:mkms2}
which can be exploited within the context of specific applications.
Meanwhile, the rate of increase under Algorithm~\ref{alg:mkms2} as a function of the number of block matrices is the smallest among all algorithms.

Figure~\ref{fig:2} presents the representative set of performance result comparisons for the case where the number $m$ of diagonal block matrices $\bP_{ii}$, $i\in[m]$, is fixed to be $20$, the magnitude $\epsilon$ of transitions in the off-diagonal block matrices $\|\bP_{ij}\| = O(\epsilon)$, $i,j\in[m], i\neq j$, is fixed to be $0.1$, and the dimensions $n_i$ of the diagonal block matrices $\bP_{ii}$, $i\in[m]$, are varied such that 
$n_i\in\{100, 200, 500, 1000, 2000\}$.
We first observe from Figure~\ref{fig:2}(a) that, once again, the mean number of iterations under Algorithm~\ref{alg:mkms} are essentially identical to that of the baseline KMS algorithm, with the mean number of iterations under Algorithm~\ref{alg:mkms2} only somewhat higher requiring a few additional iterations (though of lower computational complexity).
The mean number of iterations as a function of the dimensions $n_i$ of the diagonal block matrices $\bP_{ii}$, $i\in[m]$, remains relatively flat for all algorithms, which is as expected due to the shared properties of the outer loop of all three algorithms.
From Figure~\ref{fig:2}(b), we observe that both of 
the representative examples of
our general algorithmic
framework
provide orders of magnitude reduction in the mean computation time and that such performance improvements are consistent
and increasing
across the various dimensions $n_i$ of the diagonal block matrices $\bP_{ii}$, $i\in[m]$.
The mean computation time for each algorithm increases with $n_i$ as expected due to the linear increase in the size of the linear systems being solved, with the rate of increase somewhat larger for the baseline KMS algorithm, with somewhat smaller rates of increase for both of 
the representative examples of
our
algorithmic  
framework,
and with somewhat smaller rates of increase for Algorithm~\ref{alg:mkms2} than for Algorithm~\ref{alg:mkms}.
Once again, we note that our general approach in Algorithm~\ref{alg:mkms2} provides the lowest mean computation times for the reasons described above, although the differences in mean computation times with our general approach in Algorithm~\ref{alg:mkms} are relatively small but tend to grow somewhat with increases in the dimensions $n_i$ of the diagonal block matrices $\bP_{ii}$, $i\in[m]$,
with some differences from this trend due to the nature of the block-Jacobi preconditioner as described below.

Figure~\ref{fig:3} presents the representative set of performance result comparisons for the case where the dimensions $n_i$ of the diagonal block matrices $\bP_{ii}$, $i\in[m]$, are fixed to be $500$, the number $m$ of diagonal block matrices $\bP_{ii}$, $i\in[m]$, is fixed to be $20$, and the magnitude $\epsilon$ of transitions in the off-diagonal block matrices $\|\bP_{ij}\| = O(\epsilon)$, $i,j\in[m], i\neq j$, is varied such that $\epsilon\in\{0.01, 0.05, 0.1, 0.15, 0.2\}$.
We first observe from Figure~\ref{fig:3}(a) that the mean number of iterations under Algorithm~\ref{alg:mkms} are 
either identical to that of the baseline KMS algorithm or differ by at most one,
with the mean number of iterations under Algorithm~\ref{alg:mkms2} somewhat higher requiring several additional iterations (though of lower computational complexity).
The mean number of iterations as a function of the magnitude $\epsilon$ of transitions in the off-diagonal block matrices $\|\bP_{ij}\| = O(\epsilon)$, $i,j\in[m], i\neq j$, remains relatively flat for both KMS and Algorithm~\ref{alg:mkms}, as expected again due to the shared properties of the outer loop of all three algorithms.
However, as observed in Figure~\ref{fig:3}(a), the mean number of iterations under Algorithm~\ref{alg:mkms2} increases and requires more outer-loop iterations as $\epsilon$ increases.
Due to the nature of the block-Jacobi preconditioner, the convergence of the 
RI
method becomes slower with increasing $\epsilon$ because now the part of the iteration matrix that is left outside the block-diagonal structure becomes more computationally significant. 
From Figure~\ref{fig:3}(b), we observe that 
both of 
the representative examples of
our general algorithmic
framework
provide orders of magnitude reduction in the mean computation time with such performance improvements remaining relatively flat across the various magnitudes $\epsilon$ of transitions in the off-diagonal block matrices $\|\bP_{ij}\| = O(\epsilon)$, $i,j\in[m], i\neq j$.
The mean computation time for each algorithm remains relatively flat as a function of the magnitude $\epsilon$, which is as expected due to the shared properties of the outer loop of all three algorithms.
We note that our general approach in Algorithm~\ref{alg:mkms} provides somewhat lower mean computations times than Algorithm~\ref{alg:mkms2} for 
larger values of $\epsilon$.
This is due to the nature of the block-Jacobi preconditioner, as described above.

\begin{figure}[htbp]
\centering
\begin{subfigure}{0.37\textwidth}
    \includegraphics[width=\textwidth]{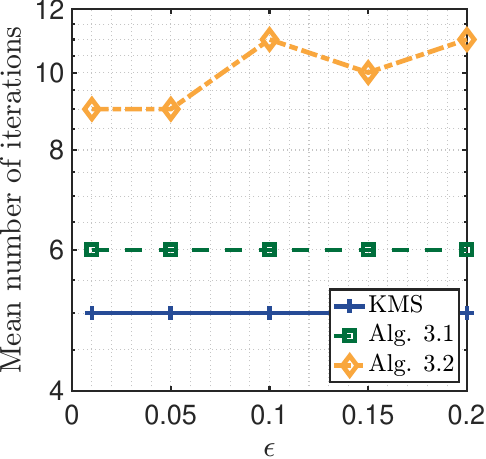}
    \caption{Mean number of iterations.}
    \label{figR:a}
\end{subfigure}
\hfill
\begin{subfigure}{0.37\textwidth}
    \includegraphics[width=\textwidth]{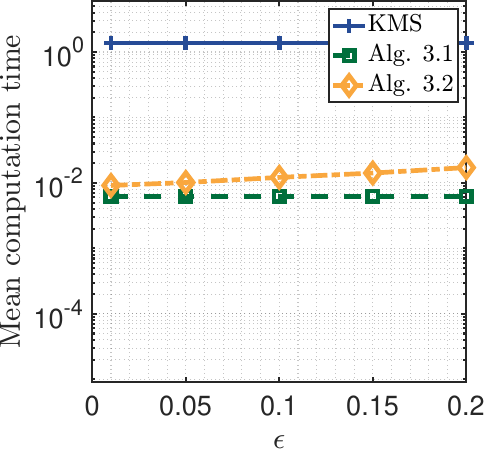}
    \caption{Mean computation times (seconds).}
    \label{figR:b}
\end{subfigure}
\caption{Performance comparison of our Algorithm~\ref{alg:mkms} and Algorithm~\ref{alg:mkms2} with the KMS baseline for diagonal matrices $\bP_{ii}$ from real-world linear systems having dimensions $n_i=730$ and number of blocks $m=20$, as a function of $\epsilon\in\{0.01, 0.05, 0.1, 0.15, 0.2\}$ for the off-diagonal matrices $\|\bP_{ij}\| = O(\epsilon)$.
The comparative performance results represent averages over $10$ independent trial runs executed on a simulated H$100$ GPU architecture.}
\label{fig:R}
\end{figure}

Next, we turn to present in Figure~\ref{fig:R} the results from a representative sample of our second collection of numerical experiments
based on selecting the $n_i\times n_i$ diagonal block matrices $\bP_{ii}$, $i\in[m]$, from a standard library of matrices for real-world linear systems~\cite{SpareMatrixLibrary}.
In particular, Figure~\ref{fig:R} presents the representative set of performance result comparisons for the case where the diagonal block matrices $\bP_{ii}$, $i\in[m]$, consist of a particular real-world matrix in computational chemistry selected from the standard library with $n_i = 730$, the number $m$ of diagonal block matrices $\bP_{ii}$, $i\in[m]$, is fixed to be $20$, and the magnitude $\epsilon$ of transitions in the off-diagonal block matrices $\|\bP_{ij}\| = O(\epsilon)$, $i,j\in[m], i\neq j$, is varied such that $\epsilon\in\{0.01, 0.05, 0.1, 0.15, 0.2\}$.
We first observe from Figure~\ref{fig:R}(a) that the mean number of iterations under Algorithm~\ref{alg:mkms} 
differs from the baseline KMS algorithm by one,
with the mean number of iterations under Algorithm~\ref{alg:mkms2} somewhat higher requiring several additional iterations (though of lower computational complexity).
The mean number of iterations as a function of the magnitude $\epsilon$ of transitions in the off-diagonal block matrices $\|\bP_{ij}\| = O(\epsilon)$, $i,j\in[m], i\neq j$, remains
flat for both KMS and Algorithm~\ref{alg:mkms}, as expected again due to the shared properties of the outer loop of all three algorithms.
However, as observed in Figure~\ref{fig:R}(a), the mean number of iterations under Algorithm~\ref{alg:mkms2} increases and requires more outer-loop iterations as $\epsilon$ increases.
Due to the nature of the block-Jacobi preconditioner, the convergence of the 
RI
method becomes slower with increasing $\epsilon$ because now the part of the iteration matrix that is left outside the block-diagonal structure becomes more computationally significant. 
From Figure~\ref{fig:R}(b), we observe that 
both of 
the representative examples of
our general algorithmic
framework
provide orders of magnitude reduction in the mean computation time with such performance improvements remaining 
relatively
flat across the various magnitudes $\epsilon$ of transitions in the off-diagonal block matrices $\|\bP_{ij}\| = O(\epsilon)$, $i,j\in[m], i\neq j$,
with some differences between Algorithm~\ref{alg:mkms} and Algorithm~\ref{alg:mkms2} due to the block-Jacobi preconditioner reasons noted above.
The mean computation times for each algorithm remains relatively flat as a function of the magnitude $\epsilon$, which is as expected due to the shared properties of the outer loop of all three algorithms.
We note that our general approach in Algorithm~\ref{alg:mkms} provides somewhat lower mean computations times than Algorithm~\ref{alg:mkms2} for 
larger values of $\epsilon$.
This is 
again
due to the nature of the block-Jacobi preconditioner, as described above. In particular, increasing $\epsilon$ makes the block off-diagonal couplings that are not captured by the preconditioner more significant.

\section{Conclusions}
We consider in this paper the problem of computing the stationary distribution of NCD Markov processes, a well-established area in the classical theory of Markov processes with broad applications in the design, modeling, analysis and optimization of computer systems and applications.
Our focus is on the design and analysis of numerical methods that address the performance bottlenecks of the most-efficient existing approaches and provide significant computational improvements to support the very large scale of current and emerging computer systems and applications as well as their potential deployment in online or real-time environments.
In particular,
we devise
a general mathematical framework of numerical solution methods
that exploits forms of advances in mixed-precision computation to significantly reduce 
the performance bottlenecks associated with solving systems of linear equations
and that exploits forms of advances in iterative approximate 
computing approaches
to 
mitigate
the impact of inaccurate computations, further reduce computation times, and guarantee convergence.
We then derive a mathematical analysis 
that rigorously establishes theoretical 
properties of our general algorithmic framework including
results 
on approximation errors, convergence behaviors, and other algorithmic 
characteristics.
Numerical experiments 
based on validated simulation of mixed-precision computation
empirically demonstrate that 
representative examples of
our general algorithmic 
framework
provide orders of magnitude improvements in the computation times to determine the stationary distribution of NCD Markov processes over the most-efficient existing numerical methods which are based on aggregation-disaggregation.

\newpage
\appendix

\renewcommand{\thesection}{\Alph{section}} 
\makeatletter
\renewcommand\@seccntformat[1]{\appendixname\ \csname the#1\endcsname.\hspace{0.5em}}
\makeatother

\section{Full Proofs of Some Theoretical Results}
\label{app:proofs}
We present in this appendix complete proofs of some of the theoretical results whose corresponding proofs in Section~\ref{sec:proofs} excluded some technical details and steps.

\subsection{Proof of Lemma~\ref{lem:4.2+4.3}}
\label{app:lem:4.2+4.3}
\begin{proof}
    From the specific supposition $\by^{(t-1)}_i = \by_i + O(\varphi)$ together with $0< \tau \leq \|\by_i\|_1 < 1$, we obtain
    \begin{equation*}
        \| \by^{(t-1)}_i \|_1 = \| \by_i \|_1 + O(\varphi)
        \qquad \mbox{ and } \qquad
        \frac{\| \by^{(t-1)}_i \|_1}{\| \by_i \|_1} = \frac{\| \by_i \|_1}{\| \by_i \|_1} + \frac{O(\varphi)}{\| \by_i \|_1} = 1 + O(\varphi) ,
    \end{equation*}
    from which it follows that
    \begin{equation*}
        \frac{\by^{(t-1)}_i}{\| \by^{(t-1)}_i \|_1} = \frac{\by_i}{\| \by_i \|_1} + O(\varphi) , \qquad i \in [m] .
    \end{equation*}
    This together with the definition of $\bQ^{(t-1)}_{ij} := ({\by^{(t-1)}_i}/{\| \by^{(t-1)}_i \|_1}) \bP_{ij} \ones$ from our Algorithm~\ref{alg:mkms} yields
    \begin{equation*}
        \bQ^{(t-1)}_{ij} = 
        \frac{\by_i}{\| \by_i \|_1} \bP_{ij} \ones + O(\varphi) \bP_{ij} \ones , \qquad \forall i,j \in [m], \; i \neq j.
    \end{equation*}
    Then, from the definition of $\bQ_{ij} := ({\by_i}/{\| \by_i \|_1}) \bP_{ij} \ones$ in the rightmost equation of~\eqref{eq:2.6} and since $\bP_{ij} = O(\epsilon)$ from~\eqref{eq:2.2} for $i \neq j$, we have
    \begin{equation*}
        \bQ^{(t-1)}_{ij} = \bQ_{ij} + O(\varphi \epsilon) , \qquad \forall i,j \in [m], \; i \neq j ,
    \end{equation*}
    which renders the desired result for $\bQ^{(t-1)}$ as an $O(\varphi \epsilon)$ approximation of $\bQ$ when $i \neq j$.
    Given that $\bQ^{(t-1)}$ and $\bQ$ are both stochastic matrices, we further conclude
    \begin{equation*}
        \bQ^{(t-1)}_{ii} = 1 - \sum_{j \in [m] : j \neq i} \bQ^{(t-1)}_{i j} = 1 - \sum_{j \in [m] : j \neq i} \bQ_{i j} + O(\varphi \epsilon) = \bQ_{ii} + O(\varphi \epsilon) , \qquad \forall i \in [m] ,
    \end{equation*}
    thus completing the proof of $\bQ^{(t-1)}$ as an $O(\varphi \epsilon)$ approximation of $\bQ$.

    Next, from the middle equation in~\eqref{eq:4.4}, we have
    \begin{equation}
        (\bI - \bQ) = \bW \begin{pmatrix} 0 & 0 \\ 0 & \bI - \bJ \end{pmatrix} \bW^{-1} .
    \tag{\ref{eq:4.13}}
    \end{equation}
    Given the theorem hypotheses $\bw^{(t-1)} = \bw^{(t-1)} \bQ^{(t-1)}$ and $\bw = \bw \bQ$, 
    where the latter follows from~\eqref{eq:2.8},
    we take the difference between the former
    equation
    and the latter
    equation
    on both sides
    and derive    
    \begin{align}
        \bw^{(t-1)} - \bw & = \bw^{(t-1)} \bQ^{(t-1)} - \bw \bQ \nonumber \\
        \bw^{(t-1)} - \bw & = \bw^{(t-1)} \bQ^{(t-1)} - \bw^{(t-1)} \bQ + \bw^{(t-1)} \bQ - \bw \bQ \nonumber \\
        (\bw^{(t-1)} - \bw) (\bI - \bQ) & = \bw^{(t-1)} ( \bQ^{(t-1)} - \bQ ) \tag{\ref{eq:4.14}} \\
        (\bw^{(t-1)} - \bw) \bW \begin{pmatrix} 0 & 0 \\ 0 & \bI - \bJ \end{pmatrix} & \stackrel{(a)}{=} \bw^{(t-1)} ( \bQ^{(t-1)} - \bQ ) \bW \nonumber,
    \end{align}
    where
    (a)
    follows upon substituting \eqref{eq:4.13} into \eqref{eq:4.14} and multiplying both sides on the right by $\bW$.
    Given the theorem
    hypothesis
    $\| \bw \|_1 = \| \bw^{(t-1)} \|_1 = 1$ (according to the rightmost equation in~\eqref{eq:2.8}) and $\| \bW \| = O(1)$ from the rightmost equation in~\eqref{eq:4.4}, together with $\bQ^{(t-1)} = \bQ_{ij} + O(\varphi \epsilon)$ established above, we obtain
    \begin{equation}
        (\bw^{(t-1)} - \bw) \bW \begin{pmatrix} 0 & 0 \\ 0 & \bI - \bJ \end{pmatrix} = \bQ^{(t-1)} - \bQ = O(\varphi \epsilon) .
    \tag{\ref{eq:4.17}}
    \end{equation}
    Letting $(\hat{\ew}_1, \hat{\ew}_2, \ldots, \hat{\ew}_m) = (\bw^{(t-1)} - \bw) \bW$, we have
    \begin{align}
        \hat{\ew}_1 = ((\bw^{(t-1)} - \bw) \bW)_1 = (\bw^{(t-1)} - \bw)\ones & = \| \bw^{(t-1)} \|_1 - \| \bw \|_1 = 0 , \tag{\ref{eq:4.18}} \\
        (\hat{\ew}_2, \ldots, \hat{\ew}_m) (\bI - \bJ) & \stackrel{(a)}{=} O(\varphi \epsilon) , \nonumber
    \end{align}
    where
    (a)
    follows from \eqref{eq:4.17}.
    Since $(\bI - \bJ)$ is nonsingular from~\eqref{eq:4.5}, we first multiply both sides of the last equation on the right by $(\bI - \bJ)^{-1}$ and then take the norm on both sides to conclude
    \begin{equation*}
        \| (\hat{\ew}_2, \ldots, \hat{\ew}_m) \| \leq O(\varphi \epsilon) \| (\bI - \bJ)^{-1} \| .
    \end{equation*} 
    From~\eqref{eq:4.5},
    we obtain
    \begin{equation}
        (\hat{\ew}_2, \ldots, \hat{\ew}_m) = O(\varphi) ,
    \tag{\ref{eq:4.20}}
    \end{equation}
    and upon combining \eqref{eq:4.18} and \eqref{eq:4.20}, we have
    \begin{equation*}
        (\bw^{(t-1)} - \bw) \bW = O(\varphi) .
    \end{equation*}
    It therefore follows from the last equation in~\eqref{eq:4.4} that $(\bw^{(t-1)} - \bw) = O(\varphi)$,
    thus completing the proof of~\eqref{eq:w-approx}.

Lastly, from~\eqref{eq:2.4}, we have
\begin{equation*}
    \by_i = \sum_{j=1}^m \by_j \bP_{ji}
    \qquad \mbox{ and } \qquad
    \by_i (\bI - \bP_{ii}) = \sum_{j=1:j\neq i}^m \by_j \bP_{ji} .
\end{equation*}
Similarly, from~\eqref{eq:4.1} and~\eqref{eq:4.1:cond3}, we obtain
\begin{equation}
    \bu_i \bP_{ii} = (1 - \et_i \epsilon + o(\epsilon) ) \bu_i
    \qquad \mbox{ and } \qquad
     \bu_i ( \bI - \bP_{ii} ) = (\et_i \epsilon + o(\epsilon) ) \bu_i ,
\tag{\ref{eq:cancel}}
\end{equation}
and therefore we conclude
\begin{align}
    (\by_i - \bu_i) (\bI - \bP_{ii}) & = \sum_{j=1:j\neq i}^m \by_j \bP_{ji} + ( -\et_i \epsilon + o(\epsilon) ) \bu_i \nonumber \\
    \by_i - \bu_i & = \left( \sum_{j=1:j\neq i}^m \by_j \bP_{ji} + ( -\et_i \epsilon + o(\epsilon) ) \bu_i \right) (\bI - \bP_{ii})^{-1} .
    \tag{\ref{eq:4.7}}
\end{align}
Given that
$\begin{pmatrix} \bu_i \\ \bU_i \end{pmatrix}$
is a nonsingular matrix of order $n_i$ and $\sum_{j=1:j\neq i}^m \by_j \bP_{ji}$ is a row vector of order $n_i$, 
there exists an order $n_i$ row vector
\begin{equation}
\bc_i = \sum_{j=1:j\neq i}^m \by_j \bP_{ji}\begin{pmatrix} \bu_i \\ \bU_i \end{pmatrix}^{-1} ,
\tag{\ref{eq:4.8}}
\end{equation}
from which it follows that $\bc_i = O(\epsilon)$ since the norm of $\begin{pmatrix} \bu_i \\ \bU_i \end{pmatrix}^{-1}$ is bounded by a constant independent of $\epsilon$
and $\sum_{j=1:j\neq i}^m \by_j \bP_{ji}$ is $O(\epsilon)$.
Then,
upon substitution of~\eqref{eq:4.8} into~\eqref{eq:4.7}, we obtain
\begin{align}
\by_i - \bu_i & = \bc_i \begin{pmatrix} \bu_i \\ \bU_i \end{pmatrix} (\bI - \bP_{ii})^{-1} + ( -\et_i \epsilon + o(\epsilon) ) \bu_i (\bI - \bP_{ii})^{-1} \nonumber \\
& \stackrel{(a)}{=} \bc_i \begin{pmatrix} (\et_i\epsilon + o(\epsilon))^{-1} & 0 \\ \bzero & (\bI - \bH_i)^{-1} \end{pmatrix} \begin{pmatrix} \bu_i \\ \bU_i \end{pmatrix} - \bu_i \nonumber \\
& = \bc_i \begin{pmatrix} \displaystyle \frac{\bu_i}{\et_i\epsilon + o(\epsilon)} \vspace*{0.15cm} \\ (\bI - \bH_i)^{-1} \bU_i \end{pmatrix} - \bu_i \nonumber \\
&= \frac{\ec_{i_1}}{\et_i\epsilon + o(\epsilon)} \bu_i + (\ec_{i_2}, \ldots, \ec_{i_{n_i}}) (\bI - \bH_i)^{-1} \bU_i - \bu_i ,
\tag{\ref{eq:4.8+}}
\end{align}
where (a) follows upon substituting~\eqref{eq:4.1:cond3} together with~\eqref{eq:cancel}. 
Given that $\bc_i = O(\epsilon)$, $(\bI - \bH_i)^{-1}$ is $O(1)$ from~\eqref{eq:4.2}
and $\bU_i$ is also $O(1)$ from above, it follows from~\eqref{eq:4.8+} that
\begin{align*}
    \by_i & = \frac{\ec_{i_1}}{\et_i\epsilon + o(\epsilon)} \bu_i + (\ec_{i_2}, \ldots, \ec_{i_{n_i}}) (\bI - \bH_i)^{-1} \bU_i \nonumber \\
    \by_i - \frac{\ec_{i_1}}{\et_i\epsilon + o(\epsilon)} \bu_i & = (\ec_{i_2}, \ldots, \ec_{i_{n_i}}) (\bI - \bH_i)^{-1} \bU_i \nonumber \\
    \by_i - \frac{\ec_{i_1}}{\et_i\epsilon + o(\epsilon)} \bu_i & = O(\epsilon) .
\end{align*}
Let $\eb_i = \ec_{i_1}/(\et_i\epsilon + o(\epsilon))$.
Since $\bc_i = O(\epsilon)$, this implies that $\ec_{i_1}$ is of order $\es_i \epsilon + o(\epsilon)$ for some $\es_i > 0$, and thus we have
\begin{equation*}
    \eb_i = \frac{\ec_{i_1}}{\et_i\epsilon + o(\epsilon)} = \frac{\es_i \epsilon + o(\epsilon)}{\et_i\epsilon + o(\epsilon)} = \frac{\es_i}{\et_i} + o(1) = O(1) ,
\end{equation*}
from which the desired result~\eqref{eq:4.6} follows.
\end{proof}

\subsection{Proof of Theorem~\ref{thm:lem4.1+thm4.1}}
\label{app:thm:lem4.1+thm4.1}
\begin{proof}
Let the error in the approximate solution $\by^{(t-1)}$ at iteration~$t-1$ of Algorithm~\ref{alg:mkms} be of order $\varphi$, namely $\by^{(t-1)} = \by + O(\varphi)$.
Then, from Lemma~\ref{lem:4.2+4.3}, we have
\begin{equation}
    \bQ^{(t-1)} = \bQ + O(\varphi \epsilon) \qquad \mbox{and} \qquad
    \bw^{(t-1)} = \bw + O(\varphi) .
\tag{\ref{eq:lemma4.1}}
\end{equation}
Starting with the $m$-th component of $\by^{(t)}$ for iteration~$t$ of Algorithm~\ref{alg:mkms}, we derive
\begin{align}
    \by_m^{(t)} & = \by_m^{(t)} \bP_{mm} + \sum_{j=1}^{m-1} s_j^{(t-1)} \hat{\by}_j^{(t-1)} \bP_{jm} \nonumber \\
    \by_m^{(t)} (\bI - \bP_{mm}) & \stackrel{(a)}{=} \sum_{j=1}^{m-1} (s_j + O(\varphi)) (\hat{\by}_j + O(\varphi)) \bP_{jm} \nonumber \\
    & = \sum_{j=1}^{m-1} s_j \hat{\by}_j \bP_{jm} + \sum_{j=1}^{m-1} (s_j O(\varphi) + \hat{\by}_j O(\varphi) + O(\varphi^2)) \bP_{jm} \nonumber \\
    & \stackrel{(b)}{=} \by_m (\bI - \bP_{mm}) + \sum_{j=1}^{m-1} O(\varphi) \bP_{jm} ,
    \tag{\ref{eq:4.25}}
\end{align}
where
(a) follows from substituting the rightmost equation in~\eqref{eq:lemma4.1} and from the $O(\varphi)$ approximation of $\by^{(t-1)}$,
and
(b) follows from~\eqref{eq:2.6} and~\eqref{eq:2.9} and from $s_j$ and $\hat{\by}_j$ both being $O(1)$ together with $\varphi < 1$.
Similarly to~\eqref{eq:4.8},
given that
$\begin{pmatrix} \bu_m \\ \bU_m \end{pmatrix}$
is a nonsingular matrix of order $n_m$
whose norm is bounded by a constant independent of $\epsilon$,
and given that $\sum_{j=1}^{m-1} O(\varphi) \bP_{jm}$ is a row vector of order $n_m$ whose components are $O(\varphi\epsilon)$, there exists a row vector $\bar{\bc}_m$ of order $n_m$ whose components are $O(\varphi\epsilon)$ such that
\begin{equation}
\sum_{j=1}^{m-1} O(\varphi) \bP_{jm} = \bar{\bc}_m \begin{pmatrix} \bu_m \\ \bU_m \end{pmatrix} .
\tag{\ref{eq:4.26}}
\end{equation}
Multiplying both sides of~\eqref{eq:4.25} by $(\bI - \bP_{mm})^{-1}$ and substituting~\eqref{eq:4.26}, we obtain
\begin{align}
    \by_m^{(t)} - \by_m & = \bar{\bc}_m \begin{pmatrix} \bu_m \\ \bU_m \end{pmatrix} (\bI - \bP_{mm})^{-1} \nonumber \\
    & \stackrel{(a)}{=} \bar{\bc}_m
    \begin{pmatrix} \displaystyle \frac{\bu_m}{\et_m\epsilon + o(\epsilon)} \vspace*{0.15cm} \\ (\bI - \bH_m)^{-1} \bU_m \end{pmatrix} \nonumber \\
    & = \left( \frac{\bar{\ec}_{m_1}}{\et_m\epsilon + o(\epsilon)} \right) \bu_m + (\bar{\ec}_{m_2}, \ldots, \bar{\ec}_{m_{n_m}}) (\bI - \bH_m)^{-1} \bU_m ,
    \tag{\ref{eq:4.27}}
\end{align}
where (a) follows upon substitution of~\eqref{eq:4.1:cond3}.
From~\eqref{eq:4.6}, we have
\begin{equation}
    \bu_m = \eb_m^{-1} \by_m + O(\epsilon) ,
    \tag{\ref{eq:4.28}}
\end{equation}
which together with~\eqref{eq:4.27}, \eqref{eq:4.28}, $\bar{\ec}_{m_k} = O(\varphi\epsilon)$, $\ec_{m_1} = O(\epsilon)$, $(\bI - \bH_m)^{-1} = O(1)$ and $\bU_m = O(1)$ yields
\begin{align*}
    \by_m^{(t)} - \by_m & = \left( \frac{\bar{\ec}_{m_1}}{\et_m\epsilon + o(\epsilon)} \right) \left(\frac{\et_m\epsilon + o(\epsilon)}{\ec_{m_1}} \by_m + O(\epsilon) \right) + (\bar{\ec}_{m_2}, \ldots, \bar{\ec}_{m_{n_m}}) (\bI - \bH_m)^{-1} \bU_m \\
    & = O(\varphi) \by_m + \frac{\bar{\ec}_{m_1} \cdot O(\epsilon)}{\et_m\epsilon + o(\epsilon)} + O(\varphi \epsilon) \\
    \by_m^{(t)} & = (1 + O(\varphi)) \by_m + O(\varphi \epsilon) .
\end{align*}
It then follows that
\begin{equation*}
    \hat{\by}_m^{(t)} = \frac{\by_m^{(t)}}{\| \by_m^{(t)} \|_1} = \frac{(1 + O(\varphi)) \by_m + O(\varphi \epsilon)}{\| (1 + O(\varphi)) \by_m + O(\varphi \epsilon) \|_1} = \frac{\by_m}{\| \by_m \|_1} + O(\varphi \epsilon) ,
\end{equation*}
from which we conclude
\begin{equation}
    \hat{\by}_m^{(t)} = \hat{\by}_m + O(\varphi \epsilon) .
    \tag{\ref{eq:4.30}}
\end{equation}

Now, arguing by induction with respect to the base case~\eqref{eq:4.30}, suppose 
\begin{equation}
    \hat{\by}_k^{(t)} = \hat{\by}_k + O(\varphi \epsilon)
    \label{eq:4.31-2}
\end{equation}
holds for all $k = m, m-1, \ldots, i+1$.
Then, for the $i$-th component of $\by^{(t)}$ for iteration~$t$ of Algorithm~\ref{alg:mkms}, we derive
\begin{align}
    \by_i^{(t)} & = \by_i^{(t)} \bP_{ii} + \sum_{j=1}^{i-1} s_j^{(t-1)} \hat{\by}_j^{(t-1)} \bP_{ji} + \sum_{j=i+1}^{m} \by_j^{(t)} \bP_{ji} \nonumber \\
    \by_i^{(t)} (\bI - \bP_{ii}) & \stackrel{(a)}{=} \sum_{j=1}^{i-1} (s_j + O(\varphi)) (\hat{\by}_j + O(\varphi)) \bP_{ji} + \sum_{j=i+1}^{m} (\by_j + O(\varphi\epsilon)) \bP_{ji} \nonumber \\
    & = \sum_{j=1}^{i-1} s_j \hat{\by}_j \bP_{ji} + \sum_{j=1}^{i-1} (s_j O(\varphi) + \hat{\by}_j O(\varphi) + O(\varphi^2)) \bP_{ji} + \sum_{j=i+1}^{m} \by_j \bP_{ji} + \sum_{j=i+1}^{m} O(\varphi\epsilon) \bP_{ji} \nonumber \\
    & \stackrel{(b)}{=} \sum_{j=1}^{i-1} \by_j \bP_{ji} + \sum_{j=i+1}^{m} \by_j \bP_{ji} + \sum_{j=1}^{i-1} (s_j O(\varphi) + \hat{\by}_j O(\varphi) + O(\varphi^2)) \bP_{ji} + \sum_{j=i+1}^{m} O(\varphi\epsilon) \bP_{ji} \nonumber \\
    & \stackrel{(c)}{=} \by_i (\bI - \bP_{ii}) + \sum_{j=1:j \neq i}^{m} O(\varphi) \bP_{ji} ,
    \label{eq:4.25-2}
\end{align}
where:
(a) follows from substituting the rightmost equation in~\eqref{eq:lemma4.1}, from the $O(\varphi)$ approximation of $\by^{(t-1)}$, and from substituting~\eqref{eq:4.31-2};
(b) follows from~\eqref{eq:2.6} and~\eqref{eq:2.9}; 
and
(c) follows from $s_j$ and $\hat{\by}_j$ both being $O(1)$ together with $\epsilon < \varphi < 1$.
Similarly to~\eqref{eq:4.26},
given that
$\begin{pmatrix} \bu_i \\ \bU_i \end{pmatrix}$
is a nonsingular matrix of order $n_i$
whose norm is bounded by a constant independent of $\epsilon$,
and given that $\sum_{j=1:j \neq i}^{m} O(\varphi) \bP_{ji}$ is a row vector of order $n_i$ whose components are $O(\varphi\epsilon)$, there exists a row vector $\bar{\bc}_i$ of order $n_i$ whose components are $O(\varphi\epsilon)$ such that
\begin{equation}
\sum_{j=1:j \neq i}^{m} O(\varphi) \bP_{ji} = \bar{\bc}_i \begin{pmatrix} \bu_i \\ \bU_i \end{pmatrix} .
\label{eq:4.26-2}
\end{equation}
Multiplying both sides of~\eqref{eq:4.25-2} by $(\bI - \bP_{ii})^{-1}$ and substituting~\eqref{eq:4.26-2}, we obtain
\begin{align}
    \by_i^{(t)} - \by_i & = \bar{\bc}_i \begin{pmatrix} \bu_i \\ \bU_i \end{pmatrix} (\bI - \bP_{ii})^{-1} \nonumber \\
    & \stackrel{(a)}{=} \bar{\bc}_i
    \begin{pmatrix} \displaystyle \frac{\bu_i}{\et_i\epsilon + o(\epsilon)} \vspace*{0.15cm} \\ (\bI - \bH_i)^{-1} \bU_i \end{pmatrix} \nonumber \\
    & = \left( \frac{\bar{\ec}_{i_1}}{\et_i\epsilon + o(\epsilon)} \right) \bu_i + (\bar{\ec}_{i_2}, \ldots, \bar{\ec}_{i_{n_i}}) (\bI - \bH_i)^{-1} \bU_i ,
    \label{eq:4.27-2}
\end{align}
where (a) follows upon substitution of~\eqref{eq:4.1:cond3}.
From~\eqref{eq:4.6}, we have
\begin{equation}
    \bu_i = \eb_i^{-1} \by_i + O(\epsilon) ,
    \label{eq:4.28-2}
\end{equation}
which together with~\eqref{eq:4.27-2}, \eqref{eq:4.28-2}, $\bar{\ec}_{i_k} = O(\varphi\epsilon)$, $\ec_{i_1} = O(\epsilon)$, $(\bI - \bH_i)^{-1} = O(1)$ and $\bU_i = O(1)$ yields
\begin{align*}
    \by_i^{(t)} - \by_i & = \left( \frac{\bar{\ec}_{i_1}}{\et_i\epsilon + o(\epsilon)} \right) \left(\frac{\et_i\epsilon + o(\epsilon)}{\ec_{i_1}} \by_i + O(\epsilon) \right) + (\bar{\ec}_{i_2}, \ldots, \bar{\ec}_{i_{n_i}}) (\bI - \bH_i)^{-1} \bU_i \\
    & = O(\varphi) \by_i + \frac{\bar{\ec}_{i_1} \cdot O(\epsilon)}{\et_i\epsilon + o(\epsilon)} + O(\varphi \epsilon) \\
    \by_i^{(t)} & = (1 + O(\varphi)) \by_i + O(\varphi \epsilon) .
\end{align*}
It then follows that
\begin{equation*}
    \hat{\by}_i^{(t)} = \frac{\by_i^{(t)}}{\| \by_i^{(t)} \|_1} = \frac{(1 + O(\varphi)) \by_i + O(\varphi \epsilon)}{\| (1 + O(\varphi)) \by_i + O(\varphi \epsilon) \|_1} = \frac{\by_i}{\| \by_i \|_1} + O(\varphi \epsilon) ,
\end{equation*}
from which we conclude
\begin{equation*}
    \hat{\by}_i^{(t)} = \hat{\by}_i + O(\varphi \epsilon) .
\end{equation*}
By induction, we therefore have
\begin{align*}
    \frac{\by_i^{(t)}}{\| \by_i^{(t)} \|_1} & = \frac{\by_i}{\| \by_i \|_1}+ O(\varphi \epsilon) , \qquad \forall i \in [m] .
\end{align*}
It then follows from the initial arguments in the proof of Lemma~\ref{lem:4.2+4.3} that
\begin{equation*}
    \by^{(t)}_i = \by_i + O(\varphi \epsilon) , \qquad \forall i \in [m] ,
\end{equation*}
which yields the desired result.

Next, 
under the process employed for selecting the precision of computations in \textbf{Step}~$\mathbf{5}$ and \textbf{Step}~$\mathbf{3}$ of Algorithm~\ref{alg:mkms} with respect to the condition number and the norm of the matrix $\bA$ of the linear system, we know that the assumptions of the mixed-precision analysis of 
the IR method
by Moler~\cite{Moler} are satisfied.
It then follows from the convergence results of Moler~\cite{Moler} that, under these conditions, the
IR
method is guaranteed to converge linearly with an approximation error that decreases by a factor of $O(\kappa(\bA))$ in each iteration, thus completing the proof.
\end{proof}

\end{document}